\documentclass[12pt,a4paper,reqno]{amsart}
\usepackage[T1]{fontenc}
\usepackage{lmodern}
\usepackage{amsmath,amsthm,amssymb,amsfonts,mathtools,mathrsfs}
\usepackage{bbm}
\usepackage[utf8]{inputenc}
\usepackage{tikz}
\usepackage{tikz-cd}
\tikzset{commutative diagrams/.cd,
mysymbol/.style={start anchor=center,end anchor=center,draw=none}
}
\newcommand\MySymb[2][\circlearrowleft]{%
  \arrow[mysymbol]{#2}[description]{#1}}
\usepackage{pst-node}
\usepackage{tikz-cd} 
\tikzcdset{scale cd/.style={every label/.append style={scale=#1},
    cells={nodes={scale=#1}}}}

\usetikzlibrary{positioning}
\usetikzlibrary{arrows.meta}
\usetikzlibrary{braids,knots}
\usepackage[shortlabels]{enumitem}
\setlist[itemize]{label=\textbullet}

\usepackage{listings}
\usepackage{csquotes}
\usepackage[titletoc]{appendix}
\usepackage[colorlinks,linkcolor=black,urlcolor=blue,citecolor=black]{hyperref}

\usepackage[top=3.2cm, bottom=3cm, left=3.5cm, right=3.5cm]{geometry}
\usepackage{cellspace}
\usepackage{algorithm}   

\theoremstyle{definition}
\newtheorem{definition}{Definition}
\numberwithin{definition}{subsection}

\theoremstyle{definition}

\numberwithin{Definition}{subsection}
\theoremstyle{remark}

\theoremstyle{remark}
\newtheorem{Remark}{Remark}
\numberwithin{Remark}{subsection}
\theoremstyle{remark}
\newtheorem*{remark}{Remark}

\theoremstyle{plain}
\newtheorem{Theorem}{Theorem}
\numberwithin{Theorem}{subsection}
\theoremstyle{plain}

\newcommand{\thistheoremname}{}
\newtheorem*{genericthm}{\thistheoremname}
\newenvironment{namedthm}[1]
  {\renewcommand{\thistheoremname}{#1}%
   \begin{genericthm}}
  {\end{genericthm}}

\theoremstyle{plain}

\theoremstyle{plain}
\newtheorem{Corollary}[Theorem]{Corollary}
\theoremstyle{plain}
\newtheorem{Lemma}[Theorem]{Lemma}
\theoremstyle{plain}

\theoremstyle{plain}
\newtheorem{Proposition}[Theorem]{Proposition}
\theoremstyle{plain}

\newcommand*{\N}{\mathbb{N}}
\newcommand*{\Z}{\mathbb{Z}}
\newcommand*{\Q}{\mathbb{Q}}
\newcommand*{\R}{\mathbb{R}}
\newcommand*{\Cc}{\mathbb{C}}
\newcommand*{\Hh}{\mathbb{H}}
\newcommand*{\F}{\mathbb{F}}
\newcommand*{\s}{\mathbb{S}}
\newcommand*{\A}{\mathbb{A}}
\newcommand*{\Gg}{\mathbb{G}}
\usepackage{mathtools}
\DeclareMathOperator{\Hom}{Hom}
\DeclareMathOperator{\GL}{GL}
\DeclareMathOperator{\GU}{GU}
\DeclareMathOperator{\U}{U}
\DeclareMathOperator{\SU}{SU}
\DeclareMathOperator{\SL}{SL}
\DeclareMathOperator{\GSp}{GSp}

\DeclareMathOperator{\h}{H}

\title[A $p$-adic interpolation of the Cogdell lift]{A $p$-adic interpolation of the Cogdell lift}
\author{Francesco Maria Iudica}
\date{}
\address{Beijing International Center for Mathematical Research, BICMR, PKU, 100080 Haidian District, Beijing, China}
\thanks{Beijing International Center for Mathematical Research (BICMR), 
Peking University, 100080 Haidian District, Beijing, China. \\
Email: iudica@bicmr.pku.edu.cn}
\subjclass[2020]{11F27, 11F33, 11F55}

\keywords{Theta lifts, $p$-adic modular forms, Hida families, Shimura varieties}
\email{iudica@bicmr.pku.edu.cn}
\begin{document}
\maketitle

\begin{abstract}
\noindent In this paper we obtain several results related to the $p$-adic interpolation of the classical Cogdell lift, mapping special cycles on Picard modular surfaces to elliptic modular forms. The results have a threefold nature: in the first part of the paper, we $p$-adically interpolate the adjoint Kudla lift, exploiting the previously constructed $\Lambda$-adic Kudla lift. In the second part, we construct higher weight cycles in Kuga-Sato varieties attached to Picard modular surfaces, and show modularity of the generating series of these cycles, thus obtaining a higher weight analogue of the Cogdell lift. Finally, we apply the formalism introduced by Loeffler to construct $p$-adic analytic cohomology classes of special cycles, whose generating series is proved to be a Hida family interpolating the Cogdell lifts in the weight and level variables.
\end{abstract}
\addtocontents{toc}{\protect\setcounter{tocdepth}{2}}
{\let\clearpage\relax \tableofcontents} 
\thispagestyle{empty}

\setcounter{section}{-1}
\section{Introduction}
This article deals with the $p$-adic interpolation of the classical Cogdell lift introduced in \cite{cogdell}, where the author constructs elliptic modular forms of weight $3$ and Nebentypus, whose Fourier coefficients are intersection multiplicities of special cycles on Picard modular surfaces. As a first instance of the $p$-adic variation of the above lift, in \cite{iudica} we constructed a $\Lambda$-adic version of the so-called Kudla lift. We recall that the Kudla lift is a holomorphic lift of elliptic modular forms to Picard modular forms, and hence is a classical incarnation of the more general theta correspondence for the dual reductive pair
\[
(\GU(1,1),\GU(2,1)).
\]
In the present paper, we first apply the results in \textit{loc. cit.} to study the $p$-adic variation of the adjoint Kudla lift, obtained from the same theta kernel, and which takes holomorphic Picard modular forms to elliptic modular forms.

As it is known from classical results of Kudla-Millson \cite{km} and Tong-Wang \cite{tongwang}, the adjoint Kudla lift essentially coincides with the Cogdell lift introduced in \cite{cogdell}. This can be seen as a perfect analogy with the Hilbert case, where the Hirzebruch-Zagier lift is adjoint to the Doi-Naganuma lift, as conjectured in \cite[p. 109]{hirzebruch}, and proved in \cite[p. 164]{zagier}. 

In the first part of the paper, we achieve a $p$-adic interpolation of the adjoint Kudla lift, which is compatible with the previously obtained interpolation of the Kudla lift, cf. \cite[Theorem 17]{iudica}, (see Corollary \ref{main1}):

\begin{namedthm}{Main Theorem A}
There exists an analytic endomorphism 
\[
\Psi:\s^{\emph{ord}}(\omega_{K/\Q})\to \s^{\emph{ord}}(\omega_{K/\Q})
\] 
such that for all weights $\kappa\in \mathfrak{X}_{\emph{alg}}(\Lambda_1)$, and every Hida family $\mathfrak{F}\in \s^{\emph{ord}}(\omega_{K/\Q})$,
\[
\Omega_\kappa^{-1} \Psi(\mathfrak{F})(\kappa)=\mathcal{L}^\dagger\mathcal{L}(\mathfrak{F}(\kappa)).
\]
\end{namedthm}
\noindent In the proof of this result, we use the construction of a $p$-adic $L$-function developed by \cite{eischen}.
We interpret this result as an arithmetic interpolation of the lift of Cogdell.

In the second part, we construct higher weight cycles $C_n^k$ in Kuga-Sato varieties attached to Picard modular surfaces, and show modularity of their generating series. The second main result we obtain is (see Theorem \ref{main2}):

\begin{namedthm}{Main Theorem B}
For $k\geq 0$, the theta series with coefficients in the cohomology of the Kuga-Sato variety $\mathcal{A}_k$
\[
\sum_{n\geq 0}C_n^k q^n\in \emph{H}^{4k+2}(\mathcal{A}_k,\Q)[[q]]
\]
is modular of weight $2k+3$, i.e. for $\varphi\in \emph{H}^{8k+2}(\mathcal{A}_k,\Q)$ with compact support, the $q$-expansion
\[
\sum_{n\geq 0}(C_n^k, \varphi) q^n
\]
is an elliptic modular form of weight $2k+3$.
\end{namedthm}
\noindent To prove the above theorem, we compare our cycles to the cycles with coefficients defined in \cite{little}, applying the formalism developed by Funke-Millson in \cite{funke}.
The major difference between Little's methods and ours is that our cycles are constructed geometrically so as to be algebraic cycles in the cohomology of some variety with constant coefficients.

By means of the adjointness property of the Kudla lift and the Cogdell lift, made explicit in Proposition \ref{analogue} below, we are able to relate the geometric construction of higher weight cycles in Kuga-Sato varieties to the higher weight adjoint Kudla lift, which has been interpolated in a $p$-adic family varying in weight and level in the first half of the paper, in terms of periods of Picard modular forms (see Corollary \ref{relation}). 
In the final part of the paper we construct a $p$-adic analytic family of unitary special cycles, by applying Loeffler's methods \cite{loeffler}. Furthermore, we compute the intersection multiplicities of these cycles under a suitably twisted "big" Poincaré pairing, which can be seen as a generalization of a construction of Ohta \cite[Proposition 4.1.13]{ohta}, see also \cite[Proposition 7.13]{fornea} or \cite[Lemma 1.1]{darmon} for the same construction in different settings. The $\Lambda$-adic intersection multiplicities are then packaged in a formal $q$-expansion, whose specializations at arithmetic points are proven to be modular. This provides us with a $\Lambda$-adic Cogdell theorem:
\begin{namedthm}{Main Theorem C}
Let $$\zeta_\infty\in e\h_{\emph{Iw}}^2(Q^0_G,\Z_p)$$ be a cuspidal big cohomology class that is orthogonal to the Chern class. The $\Lambda$-adic $q$-expansion
\[
\Phi_{\zeta}(z):=\sum_{n\geq 1}[\xi_{n,\infty}^{\emph{n.o.}}, \zeta_\infty] q^n\in \Lambda[[q]]
\]
interpolates the modular forms constructed by Cogdell and their higher weight analogues, i.e. for each weight $k> 0$ and level $r> 0$, 
\[
\nu_{2k,r}(\Phi_{\zeta_\infty^{\emph{n.o.}}}(z))\in S_{2k+3}(\Gamma_1(p^rD)).
\]
\end{namedthm}
See Theorem \ref{bigcog} below for more details.
Here the big classes $\xi_{n,\infty}^{\text{n.o.}}$ are constructed \textit{ad hoc} to interpolate the cycles $C_n^k$, so that the modular forms $\nu_{2k,r}(\Phi_{\zeta_\infty^{\emph{n.o.}}}(z))$ correspond, roughly speaking, to the generating series of the higher weight special cycles. We note here how the same methods might in principle be applied to the analogue problem for unitary Shimura varieties of signature $(m,1)$. We hope to come back to this question in a future work.

%Given the compatibility of the various constructions, classical and $p$-adic, it is natural to ask whether the cohomology classes obtained from the classical construction via Kuga-Sato varieties, agree with the $p$-adic construction obtained by the machinery developed by \cite{loeffler}. We hope to return to this question in a future work.

The paper is organized as follows: in \S$1$ we introduce the necessary notations and recall the definition of a Picard modular surface, both as a moduli space and as a Shimura variety. In \S$2$-$3$ we recall Cogdell's classical results on liftings of special cycles on Picard modular surfaces, and relate these with the adjoint of the Kudla lift. The first main result, Corollary \ref{main1}, on the $p$-adic interpolation of the adjoint Kudla lift, occupies \S$4$. The construction of higher weight cycles in Kuga-Sato varieties appears in \S$5$. Finally, the $p$-adic interpolation of the Cogdell lift can be found in \S$6$.
\section{The unitary Shimura variety of signature $(2,1)$}
Endow the $K$-vector space $V=K^3$ with the hermitian pairing
\begin{center}
    $(u,v)=\overline{u}^tJ v$,
\end{center}
for 
\[
J=\begin{pmatrix}& & \delta^{-1}\\
    & 1 & \\
    -\delta^{-1} & &\end{pmatrix}.
\]
We identify $V_\Cc=V\otimes_{K} \Cc$ with $\Cc^3$, letting $K$ act on it via $\sigma$. The space $V_\Cc$ is hermitian of signature $(2,1)$. %Furthermore, one can see that any hermitian space over $K$ whose signature at the infinity place is $(2,1)$ is isomorphic to $V$, possibly after rescaling the hermitian form by a positive rational number.
\begin{definition}
Let $G:=\GU(V,(,))$ be the general unitary group of $V$. For any ring $R$, define the $R$-points of the group scheme $G$ as,
\[
G(R):=\{(g,\mu(g))\in \GL_3(R\otimes\mathcal{O}_K) \times R^\times \mid \overline{g}{^t}Jg =\mu(g) J\}.
\]
\end{definition}
\noindent The map $\mu:G\to \mathbb{G}_m$ is a character of algebraic groups, and is called the similitude factor. Note that taking determinants of $\overline{g}{^t}Jg=\mu(g)J$ we see that 
\[
\mu(g)^3=\det \overline{g}\det g=|\det g|^2.
\]
We define $\U:=\ker \mu$. As $\mu$ is uniquely determined by $g$, we will often identify the couple $(g,\mu(g))$ with $g$. We let $B$ denote the Borel subgroup of all upper triangular matrices in $G$, which is also the stabilizer of the point $(1:0:0)$ in $\mathbb{P}^2(K)$. Its unipotent radical is denoted by $N$. We write elements of $N$ as $[v,u]$ with $v\in \text{Res}_{K/\Q}\Gg_{\text{a}}$ and $u\in \Gg_{\text{a}}$, where
\[
[v,u]=\begin{pmatrix}
1&\delta\overline{v}&u+\delta v\overline{v}/2\\
0&1&v\\
0&0&1
\end{pmatrix}.
\]
The center of $N$ is the group of all elements $[0,u]$ and is therefore isomorphic to $\Gg_{\text{a}}$.

The group $G_\infty=G(\R)$ acts on $\mathbb{P}^2_\Cc=\mathbb{P}(V_\Cc)$ by projective linear transformations and preserves the open subdomain $\mathfrak{X}$ of negative definite lines, which is biholomorphic to the open unit ball in $\Cc^2$. Every negative definite line is represented by a unique vector $(z,w,1)^t$ and such a vector represents a negative line if and only if
\begin{center}
    $\eta(z,w):= (z-\overline{z})/\delta-w\overline{w}> 0$.
\end{center}
For a point $\mathfrak{z}=(z,w)\in \mathfrak{X}$, we denote $\underline{\mathfrak{z}}=(z,w,1)^t$ the corresponding negative vector in $V$.
The upper half-plane $\h$ embeds in $\mathfrak{X}$ as the set of points with $w=0$.
Fix the point $\mathfrak{z}_0=(\delta/2,0)\in \mathfrak{X}$, and let $K_\infty$ be its stabilizer in $G_\infty$. The latter is compact modulo center: $K_\infty/Z(G_\infty)\cong K_\infty \cap \U(\R)$ is compact and isomorphic to $\U(2)(\R)\times \U(1)(\R)$. Since $G_\infty$ acts transitively on $\mathfrak{X}$, we may identify the latter with 
\[
G_\infty/K_\infty\cong \U(2,1)(\R)/(\U(2)(\R)\times \U(1)(\R)).
\] 
The pair $(G,\mathfrak{X})$ is a Shimura datum in the sense of Deligne, with reflex field $K$. Let $L$ be a full $\mathcal{O}_{K}$-lattice in $V$, such that $L$ is included in
\[
L^\vee:=\{v\in V\mid \text{tr}((u,v))\in \Z,\; \forall\, u\in L\},
\]
the dual of $L$. Let $G(L)_f\subset G(\mathbb{A}_f)$ the subgroup stabilizing $\hat{L}:=L\otimes \hat{\Z}$. An adelic level subgroup is an open compact subgroup $\mathcal{U}_f\subseteq G(L)_f$. We fix henceforth a prime $p$, which splits in $K$. We will suppose that $\mathcal{U}_f=\mathcal{U}^p\mathcal{U}_p$ with $\mathcal{U}^p\subseteq G(\mathbb{A}^{p,\infty})$ and $\mathcal{U}_p\subseteq G(\Z_p)$. In the following, we will mainly be interested in the case where $\mathcal{U}_p$ is hyperspecial, i.e. $\mathcal{U}_p=G(\Z_p)$, or parahoric. More precisely, fix the canonical basis $\{e_1,e_2,e_3\}$ for $\F_p^3$. Call $P_1\subset \GL_3(\Z_p)$ (resp. $P_2\subset \GL_3(\Z_p)$), the inverse images via the reduction map $\GL_3(\Z_p)\to \GL_3(\F_p)$ of the subgroups of $\GL_3(\F_p)$ that stabilize the line $\langle e_1\rangle$ (resp. the plane $\langle e_1,e_2\rangle$). 
Since $p$ splits in $K$, under the isomorphism $G(\Z_p)\cong \mathbb{G}_m(\Z_p)\times \GL_3(\Z_p)$, the parahoric level structure that we consider corresponds to the case $\mathcal{U}_p=\mathbb{G}_m(\Z_p)\times P_2$. 

\noindent For each choice of adelic level subgroup $\mathcal{U}_f\subseteq G(L)_f$, define $\Gamma_{\mathcal{U}_f}:=\mathcal{U}_f\cap G(\Q)$, the associated classical level subgroup. Call the conductor of $\mathcal{U}_f$, or of $\Gamma_{\mathcal{U}_f}$ indistinctly, the squarefree integer which is the product of the primes $\ell$ at which $\mathcal{U}_\ell\subsetneq G(\Z_\ell)$.
The double coset space
\[
S_{\mathcal{U}_f}(\Cc):= G(\Q)\backslash G(\A_\Q)/\mathcal{U}_fK_\infty=G(\Q)\backslash (\mathfrak{X}\times G(\A_{\Q,f}))/\mathcal{U}_f,
\]
is in fact the set of complex points of a quasi-projective variety $S_{\mathcal{U}_f}/\Cc$ of dimension two, the Shimura variety associated to the reductive group $G$, of level $\mathcal{U}_f$. The character $\nu:\det/\mu:G\to T$, where $T$ is the torus $\text{Res}_{K/\Q}\Gg_m$, induces a map
\[
\nu:S_{\mathcal{U}_f}(\Cc)\to T(\A_\Q)/T(\Q)\nu(Z(\R)\mathcal{U}_f)
\]
from $S_{\mathcal{U}_f}(\Cc)$ to a finite group, whose fibers are the connected components of $S_{\mathcal{U}_f}(\Cc)$. In particular, for $\mathcal{U}_f=G(L)_f$, the Shimura surface $S$ of level $1$ has $h_K$ connected components, cf. \cite[p. 13]{goren}. The principal congruence subgroup $\mathcal{U}_f(N)$ of $G(L)_f$ is the normal subgroup obtained as the kernel of the reduction modulo $N$ map on $G(L)_f$, i.e. the subgroup operating trivially on $\hat{L}/N\hat{L}$.

The Shimura varieties obtained this way allow an interpretation as moduli spaces for polarized abelian varieties with additional endomorphism and level structures. To sketch this interpretation, let $V=K^3$, $L=\mathcal{O}_K^3, J$, be as above. Given a point $\mathfrak{z}\in \mathfrak{X}$, we obtain a splitting $V_\Cc=V_+\oplus V_-$ as a sum of a positive two-dimensional space $V_+$ and a negative one-dimensional space $V_-$. This splitting defines a canonical complex structure $j_\mathfrak{z}$ on $V_\Cc$: for $v=v_++v_-\in V_+\oplus V_-$ and $z\in \Cc$, set
\[
j_\mathfrak{z}(z)(v_++v_-):=zv_++\overline{z}v_-.
\]
Additionally, the alternating $\R$-linear form $E=\text{tr}_{\Cc/\R}\delta(\cdot,\cdot)$ has the properties $E(j_\mathfrak{z}(i)x,j_\mathfrak{z}(i)y)=E(x,y)$ and $E(j_\mathfrak{z}(i)x,x)<0$ for $x\in V_\Cc\setminus \{0\}$. If therefore $\Lambda\subset V$ is on $\mathcal{O}_K$-lattice such that the form $E$ is integral on $\Lambda$, we obtain an abelian threefold $V_\Cc/\Lambda$, with complex structure $j_\mathfrak{z}$ and polarization given by $E$, together with an $\mathcal{O}_K$-multiplication, such that the Rosati involution induces the non-trivial automorphism of $\mathcal{O}_K$.

We now associate to a pair $(\mathfrak{z},g_f)\in \mathfrak{X}\times G(\A_{\Q_f})$ the lattice $\Lambda=g_fL$ and change the hermitian form by the factor $|\mu(g_f)|_{\A_{\Q}}$. Then, our construction gives a polarized abelian threefold $A$ with additional $\mathcal{O}_K$-structure. Multiplication of $g_f$ by an element of $G(L)_f$ on the right does not change anything, and multiplication of $(\mathfrak{z}, g_f)$ by an element of $G(\Q)$ from the left induces an isomorphism of the two structures. Consequently, we have an interpretation, over $\Cc$, of the above surface of level $1$ as a moduli space. To get a corresponding interpretation for level $N$ Picard variety, we must add a level structure, i.e. a trivialization of the $N$-torsion of $A$, cf. \cite[p. 11]{goren}.
\subsection{PEL type moduli spaces}
Since $S_{\mathcal{U}_f}(\Cc)$ parametrizes abelian varieties with additional structure over $\Cc$, it is natural to give it an arithmetic structure by showing the representability of the corresponding moduli functor problem over a number ring.

Fix a prime ideal over $p$, say $\mathfrak{p}$. For each level subgroup $\mathcal{U}_f=\mathcal{U}^p\mathcal{U}_p\subset G(\A_{\Q,f})$, we define the following moduli problem that sends each connected scheme $T$ over $Spec(\mathcal{O}_{K,\mathfrak{p}})$ to the set of tuples $(A,\lambda, \eta_{\mathcal{U}_f},i)$ up to isomorphism, where:
\begin{enumerate}[(i)]
\item $A\to T$ is an abelian scheme of relative dimension $3$.
\item $i:\mathcal{O}_K\to \text{End}(A)\otimes_\Z \Z_{(p)}$ is an injective morphism of $\Z_{(p)}$-algebras.
\item $\lambda:A\to A^\vee$ is a principal polarization such that the Rosati involution induces complex conjugation on $i(\mathcal{O}_K)$
\item $\eta_{\mathcal{U}_f}$ is a $\mathcal{U}_f$-level structure.
\end{enumerate}
We explain what the level structure $\eta_{\mathcal{U}_f}$ means. For an abelian variety $A$ of dimension $3$ over $T$ satisfying the first three conditions above, and a geometric point $x\in T$, we can look at the trivialization of the Tate module
\[
\h_1(A_x,\A_{\Q,f}^p)\cong V\otimes\A_{\Q,f}^p.
\]
This is an isomorphism of symplectic modules, with the symplectic structure on the left hand side given by the Weil pairing and on the right hand side by $E$. The group $G(\A_{\Q,f}^p)$ acts on $V\otimes\A_{\Q,f}^p$ and so does $\mathcal{U}^p\subset G(\A_{\Q,f}^p)$. We are interested in two cases:
\begin{enumerate}[(i)]
\item If $\mathcal{U}_p=G(\Z_p)$, then $\eta_{\mathcal{U}_f}$ means a choice of a $\mathcal{U}^p$-orbit of the above isomorphism.
\item If $\mathcal{U}_p=\Z_p^\times\times P_2$, then $\eta_{\mathcal{U}_f}$ means a choice of a $\mathcal{U}^p$-orbit of the above isomorphism, plus a choice of a subgroup $H$ of $A[\mathfrak{p}]$ of order $p^2$.
\end{enumerate}
The following is well-known and can be found in \cite[I.4.1.11]{lan}:
\begin{Theorem}
If the level $\mathcal{U}^p$ is neat, the above moduli problem is represented by a smooth, quasi-projective scheme $S_{\mathcal{U}_f}$ over $\mathcal{O}_{K,\mathfrak{p}}$.
\end{Theorem}
Above $S_{\mathcal{U}_f}$ there is a universal abelian scheme $\mathcal{A}$ of dimension $3$ with an $\mathcal{O}_K$-action. We also have the relative differential sheaf $\omega_{\mathcal{A}}:=e^*\Omega^1_{\mathcal{A}/S_{\mathcal{U}_f}}$, where $e:S_{\mathcal{U}_f}\to \mathcal{A}$ is the identity section. This is a locally free sheaf of rank $3$ which inherits the $\mathcal{O}_K$-action
of $\mathcal{A}$. We have a decomposition $\omega_{\mathcal{A}} = \omega_{\mathfrak{p}} \oplus \omega_{\overline{\mathfrak{p}}}$.
In addition, by our assumption on the signature of $V$, the type of $\omega_{\mathcal{A}}$ as an $\mathcal{O}_{S_{\mathcal{U}_f}}$-module is also $(2, 1)$. Thus $\text{rank}_{\mathcal{O}_{S_{\mathcal{U}_f}}}
(\omega_{\mathfrak{p}}) = 2$ and $\text{rank}_{\mathcal{O}_{S_{\mathcal{U}_f}}}
(\omega_{\overline{\mathfrak{p}}}) = 1$
on all connected components of $S_{\mathcal{U}_f}$.
\section{Cogdell's Theorem}
For a fixed lattice $L\subset V$ as above, let $\Gamma\subset G(L)_f\cap G(\Q):=\Gamma(L)$ be a neat congruence subgroup. %of matrices which are in the upper unipotent modulo $p^rD$.
Let $S$ be the complex points of the Shimura variety associated to $\GU(2,1)$, and of level ${\Gamma}$. We also let $S^{\text{res}}$ be the canonical resolution of singularities of the Baily-Borel compactification of $S$, cf. \cite[\S 2]{cogdell} or \cite[p. 15]{goren}. 

First, let us briefly recall the definitions of the cycles $C_n^c$, see \cite[p. 121]{cogdell} for the details. Given a positive definite vector $v\in L$, we write $\Gamma_v$ for the stabilizer of $v$ in $\Gamma$. Let $\mathbb{B}_v\subset \mathbb{B}_2$ be the set defined by
\[
\mathbb{B}_v\cong\{u\in v^\perp\mid \langle u,u\rangle<0\}
\]
under the isomorphism $\mathbb{B}_2\cong \mathfrak{X}$, where $\mathbb{B}_2\subset \Cc^2$ is the complex $2$-ball. We have a natural map
\[
\Gamma_v\backslash \mathbb{B}_v\to \Gamma\backslash \mathbb{B}_2\cong S.
\]
We define an irreducible special cycle $C_v$ as the resolution of singularities of the canonical compactification of $\Gamma_v\backslash \mathbb{B}_v$ in $S^{\text{res}}$.
For each positive integer $n$, we also define composite cycles
\[
C_n=\sum_{v\in \Gamma\backslash L(n)} C_v,
\]
where the sum is over a set of $\Gamma$-representatives for the subset $L(n)\subset L$ of vectors of norm $n$.
We have an orthogonal decomposition, cf. \cite[p. 125]{cogdell},
\[
\h_2(S^{\text{res}}, \Q)=\text{Im}(\h_2(S,\Q)\to \h_2(S^{\text{res}}, \Q))\oplus \langle \text{cusp cycles}\rangle,
\]
where the first summand is the Poincaré dual of the compactly supported degree $2$ cohomology of $S$, and the second is spanned by the elliptic curves sitting above each cusp in the resolution of singularities of $S$. Given a cycle $C\in \h_2(S^{\text{res}}, \Q)$, we define $C^c$ to be its projection via
\[
\h_2(S^{\text{res}},\Cc)\twoheadrightarrow \text{Im}(\h_2(S,\Q)\to \h_2(S^{\text{res}}, \Q)).
\]
Thus, for each positive integer $n$, we consider the projection $C_n^c$ of $C_n$.
Additionally, we also define $C_0^c$ as the homology class dual to one-fourth of the first Chern class on $S$. The notation is justified by the discussion preceding Theorem 3.1 in \cite{hirzebruch}, where it is showed that the latter is cohomologous to a differential form with compact support, and the same argument applies for Picard modular surfaces.

We can easily extend all of Cogdell's results in \cite{cogdell} to the case of any neat congruence subgroup $\Gamma\subset \Gamma(L)$ using the results in Kudla, cf. \cite{kudlaint}. The only part that needs a little care is \cite[Lemma 6.2]{cogdell}. Indeed, consider a set of indices $\Sigma$ in bijection with the set of cusps in the Baily-Borel compactification of $S$. Let $S_\sigma$ be a cusp cycle sitting above the cusp indexed by $\sigma$ corresponding to an isotropic vector $v_1\in K^3$. We extend $v_1$ to a basis $\{v_1,v_2,v_3\}$ with respect to which the hermitian form has matrix
\[
\begin{pmatrix}& & \delta^{-1}\\
    & 1 & \\
    -\delta^{-1} & &\end{pmatrix}.
\]
Define $L_\sigma=L\cap \langle v_2\rangle$. We will perform our computations with $\Gamma'=\Gamma(\mathfrak{N})\subset \Gamma$, the maximal principal subgroup of $\Gamma$, however, we will see that our proof goes through using any principal subgroup contained in $\Gamma$. We have a covering $\pi:\Gamma'\backslash \mathbb{B}_2\to\Gamma\backslash \mathbb{B}_2$. Define $H=\overline{\Gamma}/\overline{\Gamma}'$, $H_{\sigma}=\overline{\Gamma}_{v_1}/\overline{\Gamma}'_{v_1}$, where the bar is the quotient by the center, and to distinguish cycles on $\Gamma'\backslash \mathbb{B}_2$ from cycles on $\Gamma\backslash \mathbb{B}_2$ we will add the superscript $'$.

\begin{Lemma}
\label{extension}
$(C_n,S_\sigma)=n\cdot |D|\cdot N(\mathfrak{N}_0)\cdot |H_{\sigma}|^{-1}\cdot|L(n)\cap L_\sigma|$.
\end{Lemma}

\begin{proof}
We want to prove the following identity:
\[
(C_n, S_\sigma)=|H|^{-1}\sum_{g\in H/H_{\sigma}}(C'_n, S'_{g\cdot \sigma}).
\]
First remark that by the results in \cite[p. 120]{cogdell} we know that the sets $C'_n\cap  S'_{g\cdot \sigma}$ for varying $g\in H/H_{\sigma}$ are disjoint.
Clearly, any point in one of these sets maps to some point in $C_n\cap S_\sigma$.
Fix an element $q\in C_n\cap S_\sigma$. There are exactly $|H|$ points in $\Gamma'\backslash \mathbb{B}_2$ which map to $q$. Namely, fixing $C'_n$ and $S'_{\sigma}$ as lifts of $C_n$ and $S_\sigma$ respectively and $\tilde{q}\in C'_n\cap S'_{\sigma}$, all the liftings are in $H\tilde{q}\subset H(C'_n\cap S'_{\sigma})$. But by definition, we have
\[
g\cdot C'_n=g\sum_{v\in \Gamma'\backslash L(n)} C'_v=\sum_{v\in \Gamma'\backslash L(n)} g\cdot C'_v=\sum_{v\in \Gamma'\backslash L(n)} C'_{gv}=C'_n,
\]
as multiplication by elements $g\in H=\overline{\Gamma}/\overline{\Gamma}'$ induces a permutation on the set $\Gamma'\backslash L(n)$. Indeed, the cycles $C_n$ are fixed by the action of the full $\Gamma(L)$.

Now we apply Lemma 6.2 in \cite{cogdell}, which gives
\[
(C_n, S_\sigma)=|H|^{-1}\cdot n\cdot N(\mathfrak{N}_0)\cdot|D|\sum_{g\in H/H_{\sigma}} |L(n)\cap L_{g\cdot\sigma}|.
\]
However, the sets $L(n)\cap L_{g\cdot \sigma}$ all have cardinality $|L(n)\cap L_{\sigma}|$, as multiplication by elements in $\Gamma\subset \SU(2,1)$ does not affect the norm of vectors. Hence we get
\[
\sum_{g\in H/H_{\sigma}}|L(n)\cap L_{g\cdot \sigma}|=|H|\cdot|H_{\sigma}|^{-1}\cdot|L(n)\cap L_\sigma|
\]
and the result.
\end{proof}

\begin{Remark}
The right hand side of the formula in the above lemma seems to depend on the principal subgroup we chose to work with, as $|H_{\sigma}|$ does. However, the quantity $N(\mathfrak{N}_0)\cdot |H_{\sigma}|^{-1}$ does not. Indeed, using \cite[Lemma 2.3 ]{cogdell}, we see that for any principal subgroup $\Gamma''=\Gamma(\mathfrak{M})$, we have
\[
|\overline{\Gamma}'_{v_1}/\overline{\Gamma}''_{v_1}|=[\mathfrak{N}_0:\mathfrak{M}_0].
\]
Now, using $\Gamma''$ instead of $\Gamma'$ in the above proof, we get
\[
(C_n,S_\sigma)=n\cdot |D|\cdot N(\mathfrak{M}_0)\cdot |\overline{\Gamma}_{v_1}/\overline{\Gamma}''_{v_1}|^{-1}\cdot|L(n)\cap L_\sigma|,
\]
but 
\[
|\overline{\Gamma}_{v_1}/\overline{\Gamma}''_{v_1}|=|\overline{\Gamma}_{v_1}/\overline{\Gamma}'_{v_1}|\cdot|\overline{\Gamma}'_{v_1}/\overline{\Gamma}''_{v_1}|=|H_{\sigma}|\cdot|\overline{\Gamma}'_{v_1}/\overline{\Gamma}''_{v_1}|,
\]
whence
\[
N(\mathfrak{M}_0)\cdot |\overline{\Gamma}_{v_1}/\overline{\Gamma}''_{v_1}|^{-1}=N(\mathfrak{M}_0)\cdot |H_{\sigma}|^{-1}\cdot [\mathfrak{N}_0:\mathfrak{M}_0]^{-1}=N(\mathfrak{N}_0)\cdot |H_{\sigma}|^{-1}.
\]
\end{Remark}
For simplicity and to work in the same setting as Cogdell, for the next results in this section we suppose that $L=\mathcal{O}_K^3$. In the general case, the modular forms we obtain have level depending on the discriminant of $L$, see for example \cite[Theorem 6.4.1]{little}.
\begin{Theorem}[Cogdell]
For every positive vector $v\in L$, the function
\[
\Phi_v(z)=\frac{1}{2}e(\Gamma_v\backslash \mathbb{B}_v)+\sum_{n=1}^{\infty}(C_v, C_n^c)\;q^n,\;\;\;\;\; q=e^{2\pi iz},
\]
defines a modular form on $\Hh$ of weight $3$ and character $\omega_{K/\Q}$ for the congruence subgroup $\Gamma_0(D)$.
\end{Theorem}
While Cogdell's result only deals with a lift of an irreducible cycle $C_v$, which, roughly speaking, is an embedded modular or Shimura curve in the Picard modular surface, we shall extend this to the full $\h_2(S, \Cc)$, obtaining a result analogous to \cite[Theorem 3.1]{hirzebruch}:
\begin{Theorem}
\label{main}
For each cohomology class $C$ in the subspace of $\h_2(S^{\text{res}},\Q)$ generated by the $C_n^c, n\in \N$, the function
\[
\Phi_C(z):=\sum_{n\geq 0}(C,C_n^c)q^n
\]
is a modular form of weight $3$, level $\Gamma_0(D)$ and Nebentypus $\omega_{K/\Q}$. It is non-cuspidal for $C=C_0^c$, and a cusp form for $C$ orthogonal to $C_0^c$. The map $C\mapsto \Phi_C$ from the space generated by the $C_n^c$ to $S_3(\Gamma_0(D),\omega_{K/\Q})$ is injective.
\end{Theorem}

\begin{proof}
The proof can be carried out exactly in the same way as in \cite[p. 107]{hirzebruch} because the same general arguments apply, however we recall here the key ideas:
first of all remark that this result is in accordance with the results in \cite{cogdell} because $(C,C_0^c)=\frac{1}{2}e(C)$, where $e(C)$ is the Euler number of $C$, cf. \cite[p. 195]{kudlaint}. Extending by linearity the result of Cogdell on the space generated by the $C_n^c$ for $n>0$, we see that we are only left to show that $\Phi_{C_0^c}$ is a modular form. However this is an immediate consequence of \cite[p. 124]{km}, since the Picard surface has finite volume and $C_0^c$ has compact support.

It is also clear that if we take $C\in C_0^\perp$, then $\Phi_C$ will be a cusp form, as its constant term is $(C,C_0^c)=0$. On the contrary, $\Phi_{C_0^c}$ is never cuspidal, as the constant term will be a multiple of the volume of the Picard modular surface. Finally, the map $C\mapsto \Phi_C$ is injective by the Hodge index theorem, cf. \cite[p. 107]{hirzebruch}.
\end{proof}

\begin{Remark}
The injectivity of the above mapping is equivalent to the fact that $\Phi_C$ is well defined for any $C\in \h_2(S^{\text{res}},\Cc)$, as this is simply zero for $C$ belonging to the orthogonal complement of the subspace generated by the special cycles $C_n^c, n\in \N$.
\end{Remark}

\section{The adjoint Kudla lift}
We turn to the study of the adjoint lift $\mathcal{L}^\dagger$ of $\mathcal{L}$, constructed from the same theta kernel $\Theta$, i.e.
\[
\mathcal{L}^\dagger(F)=(F,\overline{\Theta}), \;\;\; \forall F\in \mathscr{S}(k,\Gamma),
\]
whereas
\[
\mathcal{L}(f)=(f,\Theta), \;\;\; \forall f\in S_{k-1}(\Gamma^0,\omega_{K/\Q}).
\] 
In particular, we describe the composition $\mathcal{L}^\dagger\circ \mathcal{L}$. For a fixed even weight $k\geq 4$ and $\Gamma^0\subset \SL_2(\Z)$ a congruence subgroup, we write $\mathcal{L}_k:S_{k-1}(\Gamma^0,\omega_{K/\Q})\to \mathscr{S}(k,\Gamma)$, for the Kudla lift.

We wish to relate $\mathcal{L}_4^\dagger$ to the geometric lift
\begin{align*}
\h_{1,1}(S, \Cc)&\to M_3(\Gamma_0(D),\omega_{K/\Q})\\
[C]&\mapsto\sum_{n\geq 0}(C,C_n^c)q^n,
\end{align*}
which has been studied in \cite{cogdell}, by means of the theory developed in \cite{km} and in \cite{tongwang}.

Technically speaking, the definitions of $\mathcal{L}$ and $\mathcal{L}^\dagger$ break down for weight $k<6$, for convergence reasons, but in \cite{tongwang} it is showed that by analytic continuation the same lifts can be defined for $k=4$, corresponding to weight $3$ elliptic cusp forms. 

On the ball model of the symmetric space associated to $\GU(2,1)$, the K\"ahler metric $\mu_{Hyp}$ can be expressed by the following formula:
\[
\mu_{Hyp}(z):=2\frac{1}{(1-|z|^2)^2}.
\]
We let $d\mu_{Hyp}^{\text{vol}}$ denote the associated volume form, cf. \cite[p. 4]{baskar}. In \emph{loc. cit.} the Petersson inner product of two Picard cusp forms $F,G\in \mathscr{S}(k,\Gamma)$, is defined by the formula:
\[
( F, G):=\int_{\mathcal{F}_{\Gamma}}(1-|z|^2)^kF(z)\overline{G(z)}d\mu_{Hyp}^{\text{vol}}(z),
\]
where $\mathcal{F}_{\Gamma}$ denotes a fundamental domain of $\Gamma\backslash \mathbb{B}_2$.

\subsection{The analogue of a result of Zagier and Oda}
For $F\in \mathscr{S}(k,\Gamma)$, we want to compute the adjoint Kudla lift of $F$:
\[
\mathcal{L}^\dagger(F)=(F,\overline{\Theta}),
\]
where $\Theta$ is the theta kernel of the Kudla lift used in \cite[p. 13]{iudica}.
\begin{Proposition}
\label{analogue}
We have
\[
\mathcal{L}^\dagger(F)(z)=c_k(L)\sum_{v\in L^+}\Big(\int_{C_v}F(\mathfrak{z})(v,\underline{\mathfrak{z}})^{k}d\mu_{v}(\mathfrak{z})\Big)q^{(v,v)},
\]
where $d\mu_v$ denotes the pull-back of the measure $d\mu_{Hyp}^{\text{vol}}$ along $\mathfrak{X}_v\to \mathfrak{X}$, and $c_k(L)\in \Cc$ is a constant independent of $F$.
\end{Proposition}

\begin{proof}
Write
\[
\mathcal{L}^\dagger(F)(z)=\int_{\Gamma\backslash G}\Theta(z,g)\varphi(g)dg,\quad \text{with   } \varphi(g):=j(g,\mathfrak{z}_0)^{-k}F(g\mathfrak{z}_0).
\]
Then 
\[
\mathcal{L}^\dagger(F)(z)=\int_{\Gamma\backslash G} y^2\sum_{v\in L}f(z,g^{-1}v)\varphi(g)dg.
\]
Since $F$ is a cusp form, this integral is absolutely convergent, and we can interchange integration and summation. For each $v\in V$, there exists a Haar measure $d\mu_v^1(h)$ on $G_v$, and a measure $d\mu_v^2(g)$ on $G_v\backslash G$, satisfying
\[
\int_G\psi(g)dg=\int_{G_v\backslash G}\Big(\int_{G_v}\psi(hg)d\mu_v^1(h)\Big)d\mu_v^2(g),
\]
for any integrable function $\psi$ on $G$, see \cite[p. 108]{oda}. Then we have
\[
\mathcal{L}^\dagger(F)(z)=\sum_{v\in L}\int_{G_v\backslash G}y^2f(z,g^{-1}v)\Big(\int_{\Gamma_v\backslash G_v} \varphi(hg)d\mu_v^1(h)\Big)d\mu_v^2(g)
\]

We turn then to the evaluation of the inner integrals. First of all, the terms corresponding to $v$ with $(v,v)\leq 0$ vanish because $\mathcal{L}^\dagger(F)$ is holomorphic and cuspidal, see \cite[p. 13]{kudlaq}.
Define 
\[
b(g,v)=(v,\underline{g\mathfrak{z}_0}).
\]
Then, for any element $h\in G_v$, we have the equality
\[
j(hg,\mathfrak{z}_0)^{-k}=\frac{b(hg,v)^k}{b(1,g^{-1}v)^k}.
\]
Hence
\[
\int_{\Gamma_v\backslash G_v}\varphi(hg)d\mu_v^1(h)=b(1,v^tg^{-1})^{-k}\int_{\Gamma_v\backslash G_v}F(hg\mathfrak{z}_0)b(hg\mathfrak{z}_0,v)^kd\mu_v^1(h).
\]
Put $\mathfrak{z}=g\mathfrak{z}_0$ and let $\psi(\mathfrak{z})$ be the integral
\[
\int_{\Gamma_v\backslash G_v}F(h\mathfrak{z})b(h\mathfrak{z},v)^kd\mu_v^1(h).
\]
Then $\psi$ is a holomorphic $G_v$-invariant function on $\mathfrak{X}$, hence constant.
Thus we are left with
\[
\psi(\mathfrak{z}_0)=\int_{\Gamma_v\backslash G_v}F(h\mathfrak{z}_0)b(h\mathfrak{z}_0,v)^kd\mu_v^1(h)=\int_{\Gamma_v\backslash G_v/K_v}F(h\mathfrak{z}_0)b(h\mathfrak{z}_0,v)^kd\nu_v(h),
\]
where $K_v$ is $K\cap G_v$, and $d\nu_v$ is a measure on $\mathfrak{X}_v=G_v/K_v$ invariant under $G_v$, satisfying $d\mu_v^1(h)=d\nu_vd\kappa_v$, for the Haar measure $d\kappa_v$ on $K_v$ such that
\[
\int_{K_v}d\kappa_v=1.
\]
By uniqueness of the Haar measure, $d\kappa_v$ coincides with the pull-back of $d\mu_{Hyp}$ along $\mathfrak{X}_v\to \mathfrak{X}$, up to a constant multiple. The rest of the computation follows easily from the proof in \cite[p. 112]{oda}.

\end{proof}

Now, for $v\in L$ with $(v,v)<0$ and $\mathfrak{z}=(z,u)\in \mathfrak{X}$, let
\[
b_v(\mathfrak{z}):=(v,\underline{\mathfrak{z}})^{-1}=(\delta^{-1}\overline{v_3}z+\overline{v_2}u-\delta^{-1}\overline{v_1})^{-1}.
\]
By construction, $b_v(\mathfrak{z})$ is holomorphic in $\mathfrak{z}$, and for $g\in G(\R)$ we have
\[
b_{gv}(\mathfrak{z})=j(g^\iota,\mathfrak{z})^{-1}b_v(g^\iota\mathfrak{z}),
\]
where we define $g^\iota$ by
\[
gg^\iota=\mu(g).
\]
From this it easily follows that, for $\Gamma\subset G(\Q)$ a congruence subgroup,
\[
\theta_v(\mathfrak{z}):=\sum_{\gamma\in \Gamma/\Gamma_v}b_{\gamma v}(\mathfrak{z})^{k}\in \mathscr{S}(k,\Gamma).
\]
and
\[
\theta_n(\mathfrak{z}):=\sum_{v\in \Gamma\backslash L(-n)}\theta_v(\mathfrak{z})\in \mathscr{S}(k,\Gamma).
\]

\noindent For $n>0$ an integer, let $\phi^k_n:=\phi_n$ denote the degree $n$ Poincaré series in $S_k(\Gamma_0(DN),\omega_{K/\Q})$.

\begin{Proposition}[\text{\cite[p. 13]{kudlaq}}]
\[
\mathcal{L}(\phi_n)(\mathfrak{z})=(4\pi n)^{-k}(k-1)!\,\theta_n(\mathfrak{z}).
\]
\end{Proposition}
\noindent By a well-known property of the Poincaré series, we have
\begin{align*}
\frac{(k-2)!}{(4\pi n)^{k-1}}c_k(L)&\sum_{(v,v)=n}\int_{C_v}\Big(F(\mathfrak{z})(v,\underline{\mathfrak{z}})^{k}d\mu_{v}(\mathfrak{z})\Big)\\&=(\mathcal{L}^\dagger(F), \phi_n)=(F,\mathcal{L}(\phi_n))=\frac{(k-1)!}{(4\pi n)^{k}}(F,\theta_n),
\end{align*}
which implies
\[
c_k(L)\sum_{(v,v)=n}\int_{C_v}\Big(F(\mathfrak{z})(v,\underline{\mathfrak{z}})^{k}d\mu_{v}(\mathfrak{z})\Big)=\frac{k-1}{4\pi n}(F,\theta_n).
\]
As a first interpretation of these equalities, we see that the adjoint Kudla lift coincides with the Cogdell lift, in the sense that the following diagram commutes:

\begin{center}
\begin{tikzpicture}
    % set up the nodes
    \node (E) at (0,0) {$S_{3}(\Gamma^0)$};
    \node[right=2cm of E] (F) {$\mathscr{S}(4,\Gamma)$};
    \node[below= of F] (G) {$\textbf{SC}^\vee$};
    % draw arrows and text between them
    \draw[->] (E)--(F) node [midway,above] {$\mathcal{L}_4$};
    %\draw[->] (E)--(G) node [midway,below left] {$\tilde{\mathcal{L}_4}$};
    \draw[->] (F)--(G) node [midway,right, pos=0.44] {$\omega$};
    \draw [->,purple] (F.north) to [out=135,in=45] node [midway,above] {$\mathcal{L}_4^\dagger$} (E.north);
     \draw [->,purple] (G)--(E) node [midway,below] {$\tilde{\mathcal{L}_4}^\dagger$};
%\draw [->,purple] (G) to [out=180,in=270] node [midway,below left] {$\tilde{\mathcal{L}_4}^\dagger$} (E.south) ; old version with 5 arrows
\end{tikzpicture}
\end{center}

\noindent where the map $\tilde{\mathcal{L}}_4$ is defined as in  \cite[p. 497]{tongwang}, and the vertical map is the association
\begin{align*}
\omega:\mathscr{S}(4,\Gamma)&\to \textbf{SC}^\vee\\
F&\mapsto (\omega_F:C_v\mapsto (v,v)^{-1}(F,\theta_v)),
\end{align*}
for $\textbf{SC}\subset \h_2(S, \Cc)$ the subspace of special cycles.

\begin{Proposition}
The above map $\omega:\mathscr{S}(4,\Gamma)\to \emph{\textbf{SC}}^\vee$ induces an embedding, still denoted with the same notation:
\[
\omega:\emph{Im}(\mathcal{L}_4)\hookrightarrow \h^2(S, \Cc).
\]
 
\end{Proposition}

\begin{proof}
This follows easily from the splitting
\[
\h_2(S, \Cc)=\textbf{SC}\oplus \textbf{SC}^\perp,
\]
which gives us an embedding
\[
\textbf{SC}^\vee\hookrightarrow \h_2(S, \Cc)^\vee\cong \h^2(S, \Cc)
\]
by extending to $0$ on $\textbf{SC}^\perp$.
\end{proof}

\begin{Corollary}
We have
\[
\mathcal{L}_4^\dagger =\frac{3}{4\pi}\tilde{\mathcal{L}}_4^\dagger\circ\omega.
\]
\end{Corollary}

\noindent Finally, we quote a result due to Tong and Wang:
\begin{Theorem}[\text{\cite[p. 497]{tongwang}}]
\label{3}
For $f\in S_3(\Gamma^0,\omega_{K/\Q})$ and $C\in \h_{1,1}(S^{\text{res}},\Cc)$ such that $C\perp C_0^c$, we have:
$$\tilde{\mathcal{L}_4}(f)=\frac{i}{2}\langle f, \Omega(\tau, z)\rangle,\;\;\;\;\tilde{\mathcal{L}_4}^\dagger(C)=\frac{i}{2}\int_{C} \Omega(\tau, z),$$
where $$\Omega(\tau,z)=\sum_{n>0} [C_n^c]^\vee q^n$$ is the generating series of the Poincaré dual classes of the special cycles.
\end{Theorem}
\begin{Remark}
We can restate Theorem \ref{3} by saying that the adjoint of the Kudla lift is the unitary Kudla-Millson lift, which sends $(1,1)$-cohomology classes on the Picard modular surface to elliptic cusp forms whose Fourier coefficients are intersection numbers. Notice that this is exactly the analogue of Conjecture 2 in \cite[p. 109]{hirzebruch}, which states the adjointness of the Doi-Naganuma and the Hirzebruch-Zagier lifts. The higher weight analogue of this conjecture is proved by Zagier in \cite[\S 6]{zagier}.
\end{Remark}
\noindent We can thus reinterpret $\mathcal{L}_4^\dagger\mathcal{L}_4$ as
\[
\mathcal{L}_4^\dagger \mathcal{L}_4(f)=\frac{3}{4\pi}\sum_{n\geq 1}(C_f,C_n^c)q^n,
\]
where $C_f$ is the class $\omega(\mathcal{L}_4(f))\in \h^2(S, \Cc)$. 

Having established the geometrical interpretation, the following result allows us to describe $\mathcal{L}^\dagger \mathcal{L}$ from a purely arithmetical point of view:
\begin{Proposition}
Suppose $f\in S_k(\Gamma^0)$ is a Hecke eigenform. Then $f$ is an eigenvector of the endomorphism $\mathcal{L}^\dagger \mathcal{L}$ of $S_k(\Gamma^0)$.
\end{Proposition}

\begin{proof}
The statement is equivalent to the fact that $\mathcal{L}^\dagger \mathcal{L}(f)$ and $f$ have the same Hecke eigenvalues. This is easily proved using the proof of Theorem 12 in \cite{iudica}.
\end{proof}

\begin{Remark}
While studying the adjoint $\mathcal{L}_4^\dagger$ of $\mathcal{L}_4$ via the composition $\mathcal{L}_4^\dagger\mathcal{L}_4$ might at first seem a bit artificial, Theorem \ref{main} together with the above proposition tell us that this is the right thing to do. Indeed, take $F\in \mathscr{S}(4,\tilde{\Gamma})$, and write $F=\mathcal{L}_4(f)+F'$, with $F'$ orthogonal to the space spanned by the Kudla lift. Now, the space spanned by the special cycles coincides with the space spanned by the theta lifts, as a consequence of the main result in \cite[p. 468]{tongwang}, cf. \cite[p. 214]{soudry}. In other words, for $F'\in \mathscr{S}(4,\tilde{\Gamma})$ $$\omega_{F'}=0\in \textbf{SC}^\vee \Leftrightarrow\mathcal{L}_4^\dagger(F')=0.$$ 
Then we see that
\[
\mathcal{L}^\dagger_4(F)=\mathcal{L}_4^\dagger\mathcal{L}_4(f),
\]
so that $\mathcal{L}_4^\dagger$ actually factors through the image of the Kudla lift.
\end{Remark}

By the above result, we know that for $f$ a Hecke eigenform, there exists an eigenvalue $\lambda_f\in \Cc$ such that
\[
\mathcal{L}^\dagger \mathcal{L}(f)=\lambda_f f.
\]
This value can be characterized as follows:
\[
\lambda_f\langle f,f\rangle=\langle \mathcal{L}^\dagger \mathcal{L}(f),f\rangle=\langle \mathcal{L}(f),\mathcal{L}(f)\rangle\Rightarrow \lambda_f=\frac{\langle \mathcal{L}(f),\mathcal{L}(f)\rangle}{\langle f,f\rangle}.
\]
The study of this ratio is a classical problem, and the so-called Rallis inner product formula identifies it with a special value of an $L$-function attached to $f$, multiplied by a finite product of bad local zeta integrals and explicit correction factors, cf. \cite[p. 226]{harris}.

To be more precise, there exists an Eisenstein series $E((g,g'),\chi,s)$, depending on the complex variable $s\in \Cc$ and on a unitary Hecke character of $K$ of weight $k$ as in \cite[p. 26]{iudica}, for the Shimura variety attached to $U(2,2)$, cf. \cite[p. 231]{harris}, for which we define
\[
Z(s,f,\chi)=\int_{H(\Q)\backslash H(\A_\Q)}E((g,g'),\chi,s)f(g)\overline{f(g')}dgdg',
\]
with $dg=dg'$ the Tamagawa measure on $U(1,1)(\A_\Q)$ and $f$ a holomorphic modular form for $U(1,1)$. Here, we have set $H:=U(1,1)\times U(1,1)\subset U(2,2)$. Moreover, we can define a partial zeta-integral $Z_S(s,f,\chi)$ on a finite set of places $S$. For our purposes, we will take 
\[
S=\{p\}\cup\{\text{primes dividing } D\},
\]
where we recall that $p$ is a fixed rational prime which splits in $K=\Q(\sqrt{-D})$.
Then we have the following basic identity, due to Piatetski-Shapiro and Rallis, cf. \cite[p. 82]{eischen}:
\begin{equation}
\label{BI}
d^S(s,\chi)Z(s,f,\chi)=\langle f,f\rangle Z_\infty(s,f,\chi) Z_S(s,f,\chi)L^S(s+\frac{1}{2},\pi_f,\chi),
\end{equation}
where 
\[
d^S(s,\chi)=L^S(2s+1,\chi)L^S(2s+2,\chi\omega_{K/\Q})
\]
and $L^S(s+\frac{1}{2},\pi_f,\chi)$ is the restricted $L$-function attached to the twist of $\pi_f$ by $\chi$, where the former is the automorphic representation associated to $f$. Note that if $f$ is a modular form for $\GL_2$ and $f_K$ is the corresponding modular form for $U(1,1)$, then, by \cite[p. 242, 245]{harris}, 
\[
L(s,f_K,\chi)=L(s,f,\chi)L(s,\overline{f},\chi^{-1}).
\]
\begin{Theorem}[\text{\cite[p. 233]{harris}}]
The special value $Z(\frac{1}{2},f,\chi)$ is an algebraic number.
\end{Theorem}
\section{The $p$-adic interpolation of the adjoint Kudla lift}
\subsection{Analytic weight spaces}
We briefly recall some notions from Hida theory.
Let $\Gamma=1+p\Z_p$ be the maximal torsion-free subgroup of $\Z_p^\times$. We fix a topological generator $\gamma\in \Gamma$ such that $\Gamma=\gamma^{\Z_p}$. Let $\Lambda=\Z_p[[\Gamma]]$ and $\Lambda_1=\Z_p[[\Z_p^\times]]$ be the completed group rings on $\Gamma$ and on $\Z_p^\times$ over $\Z_p$, respectively. Clearly $\Lambda_1$ has a natural $\Lambda$-algebra structure induced from $\Lambda_1\cong \Lambda[\mu_{p-1}]$, where $\mu_{p-1}$ denotes the maximal torsion subgroup of $\Z_p^\times$ consisting of $(p-1)$-th roots of unity.

\begin{definition}[$\Lambda$-adic weight spaces]
For each $\Lambda_1$-algebra $R$ finite flat over $\Lambda$, we define the $\Lambda$-adic weight space $\mathfrak{X}(R)$ associated with $R$ as
\[
\mathfrak{X}(R):=\text{Hom}_{\text{cont}}(R,\overline{\Q}_p),
\]
on which the following arithmetic notions are introduced:
\begin{enumerate}[(i)]
\item A point $P\in \mathfrak{X}(R)$ is said to be arithmetic if there exists an integer $k\geq 2$ such that the restriction of $P$ into $\mathfrak{X}(\Lambda)\cong \Hom_{\text{cont}}(\Gamma,\overline{\Q}_p^\times)$ corresponds to a continuous character $P_k:\Gamma\to \overline{\Q}_p^\times$ satisfying $P_k(\gamma)=\gamma^k$. We denote by $\mathfrak{X}_{\text{alg}}(R)$ the set of all arithmetic points in $\mathfrak{X}(R)$.
\item An arithmetic point $P\in \mathfrak{X}_{\text{alg}}(R)$ is said to have signature $(k,\varepsilon)$ if there exists an integer $k\geq 2$ and a finite character $\varepsilon:\Z_p^\times\to \overline{\Q}_p^\times$ such that $P$ lies over the point $P_{k,\varepsilon}\in \mathfrak{X}_{\text{alg}}(\Lambda_1)\cong \Hom_{\text{cont}}(\Z_p^\times,\overline{\Q}_p^\times)$ corresponding to the character $P_{k,\varepsilon}(y)=y^k\varepsilon(y)$ on $\Z_p^\times$. For simplicity, we just write $P=(k,\varepsilon)$ and often refer to it as the arithmetic point of weight $k$ and Nebentypus $\varepsilon\omega^{-k}$.
\end{enumerate}
\end{definition} 

\subsection{A $p$-adic Rallis inner product formula}
\noindent Let 
\[
L^p(1,\pi_f,\chi)=L^S(1,\pi_f,\chi)Z^p_S(\frac{1}{2},f,\chi)
\]
\begin{Theorem}[\text{\cite[p. 243]{harris}, \cite[Corollary 9.2.1]{eischen}}]
Given a Hida family $\mathfrak{F}$, with specializations $\{f_P\}_{P\in \mathfrak{X}_{\emph{alg}}(\Lambda_1)}$, there exists an analytic function 
\[
L_p(1,\mathfrak{F}):\mathfrak{X}(\Lambda_1)\to \Cc_p
\]
and a family of periods $\Omega:\mathfrak{X}(\Lambda_1)\to\Cc_p$, such that
\[
L_p(1,P):=L_p(1,\mathfrak{F})(P)=Z_\infty(\frac{1}{2},f_P,\chi_P)Z_p(\frac{1}{2},f_P,\chi_P)L^p(1,\pi_{f_P},\chi_P)/\Omega_P
\]
for all $P\in \mathfrak{X}_{\text{alg}}(\Lambda_1)$.
\end{Theorem}
\noindent By the work of Ichino on the Siegel-Weil formula for unitary groups, \cite{ichino}, cf. \cite[p. 246]{harris}, we have the identity:
\[
\langle\mathcal{L}(f),\mathcal{L}(f)\rangle=\int_{H(\Q)\backslash H(\A_\Q)} E((g,g',\chi,\frac{1}{2}))f(g)\overline{f(g')}dgdg'.
\]
Applying this formula to a Hida family and combining it to the Basic Identity (\ref{BI}), yields the Rallis inner product formula:
\[
d^S(\frac{1}{2},\chi_P)\langle\mathcal{L}(f_P)\mathcal{L}(f_P)\rangle=\langle f_P, f_P\rangle L_p(1,P)\Omega_P,
\]
which can be rewritten as
\[
\frac{\langle\mathcal{L}(f_P)\mathcal{L}(f_P)\rangle}{\langle f_P, f_P\rangle}=d^S(\frac{1}{2},\chi_P)^{-1}L_p(1,P)\Omega_P.
\]
This is the $\lambda_{f_P}$ introduced above.

Let $\s^{\emph{ord}}(\omega_{K/\Q})$ be the free $\Lambda_1$-module of ordinary $\Lambda$-adic forms of tame level $D$ and character $\omega_{K/\Q}$.
\begin{Corollary}
\label{main1}
There exists an analytic endomorphism 
\[
\Psi:\s^{\emph{ord}}(\omega_{K/\Q})\to \s^{\emph{ord}}(\omega_{K/\Q})
\] 
such that for all $P\in \mathfrak{X}_{\emph{alg}}(\Lambda_1)$,
\[
\Omega_P^{-1} P(\Psi(\mathfrak{F}))=\mathcal{L}^\dagger\mathcal{L}(P(\mathfrak{F})).
\]
\end{Corollary}

\begin{proof}
By \cite[Proposition 6.5]{harris}, there exists an analytic function $d^S_{\text{alg}}(\frac{1}{2},\Xi):\mathfrak{X}(\Lambda_1)\to \Cc_p$, depending on the family of Hecke characters introduced in \cite[Lemma 16]{iudica}, which interpolates the values $d^S(\frac{1}{2},\chi_P)/\Omega_P^2$, for $P\in \mathfrak{X}_{\text{alg}}(\Lambda_1)$. Given a family of eigenforms $\mathfrak{F}\in \s^{\emph{ord}}(\omega_{K/\Q})$, we define 
\[
\lambda_{\mathfrak{F}}:=d^S_{\text{alg}}\Big(\frac{1}{2},\Xi\Big)^{-1} L_p(1,\mathfrak{F})
\]
and
\[
\Psi(\mathfrak{F}):=\lambda_{\mathfrak{F}}\cdot \mathfrak{F}.
\]
For the fixed $\mathfrak{F}$, this has the wanted property, by construction.
Extending by linearity on $\s^{\emph{ord}}(\omega_{K/\Q})$, we get the result.
\end{proof}

\begin{Remark}
The above corollary shows that there exists an analytic family of modular forms which interpolates 
\[
\mathcal{L}_4^\dagger \mathcal{L}_4(f)=\frac{3}{4\pi}\sum_{n\geq 1}(C_f,C_n^c)q^n
\]
in a way which is compatible with the interpolation of  the lifts $\{\mathcal{L}_k(f)\}_k$. Indeed, a $\Lambda$-adic family of such lifts, compatible with the $p$-adic variation of $f$, has been constructed in \cite{iudica}. As by definition $C_f$ is the class $\omega(\mathcal{L}_4(f))$, we interpret this as an arithmetic $p$-adic interpolation of Cogdell's main theorem in \cite{cogdell}.
\end{Remark}

\section{Fiber varieties}
In this section, we construct higher weight cycles $C_n^k$ on Kuga-Sato varieties attached to Picard modular surfaces, show modularity of their generating series, and relate the latter to the adjoint Kudla lift. Finally, in \S$6$, these higher weight cycles will be instrumental in proving modularity of the "big" generating series.

On the hermitian space $V=K^3$ as above, we define a skew-symmetric form by $$\langle u,v \rangle:=\text{Im}(u,v).$$ This induces an embedding $$\sigma:G\hookrightarrow \GSp(6):=\GSp(V, \langle\cdot,\cdot\rangle),$$ considering $V(\R):=V\otimes_\Q \R$ as a $6$-dimensional real symplectic vector space. Let $L\subset V$ be an $\mathcal{O}_K$-lattice of rank $3$ and define, for $\Gamma\subset G(\Q)$ a torsion-free congruence subgroup preserving $L$, the quotient
\[
\mathcal{A}_{\Gamma}^{(1)}:=(\Gamma\ltimes L)\backslash(\mathfrak{X}\times V(\R)).
\]
This is a complex quasi-projective variety, which can be characterized as the complex points of the universal abelian $3$-fold associated to the moduli problem attached to the Shimura variety corresponding to $G$ and of level $\Gamma$. For a non-negative integer $r$, we let
\[
\mathcal{A}:=\mathcal{A}_{\Gamma}^{(r)}:=\underbrace{\mathcal{A}_{\Gamma}^{(1)}\times_{S_{\Gamma}}\cdots\times_{S_{\Gamma}}\mathcal{A}_{\Gamma}^{(1)}}_\text{$r$ factors}
\]
denote the $r$-fold fiber product.
Then $\mathcal{A}$, too, is a smooth complex quasi-projective variety, of complex dimension $3r+2$. As a $C^\infty$ manifold
\[
\mathcal{A}\cong (\Gamma\ltimes L^r)\backslash(\mathfrak{X}\times V^r(\R)),
\] 
where $V^r$ denotes the $r$-fold direct sum of $V$ with itself. We also allow $r=0$, for which we will mean $\mathcal{A}=S_\Gamma:=S$. 

For a point $x\in \mathfrak{X}$, we let $j(x)$ be the corresponding complex structure in the Siegel space associated to $\GSp(6)$, under the embedding $\sigma$. Then $J(x):=\sigma^r\circ j(x)$ determines a complex structure on $V^r(\R)$. Gluing together the disjoint union of $\{(V^r(\R), J(x)), x\in \mathfrak{X}\}$ gives $\mathfrak{X}\times V^r(\R)$ a holomorphic structure as a vector bundle over $\mathfrak{X}$. In particular, the natural map 
\[
\varphi: \mathcal{A}\to S
\]
is holomorphic, and thus algebraic.

We will suppose that either $L$ is integral, or $L$ is the dual of an integral lattice.
If $L$ is integral, then $\langle,\rangle$ is $\Z$-valued on $L$ and $(v_1,v_2)\mapsto \langle v_1, j(x)v_2\rangle$ is a symmetric positive-definite bilinear form. This implies that $A^r(x)=(V^r(\R)/L^r, J(x), \langle, \rangle)$ is a polarized abelian variety for each $x\in \mathfrak{X}$, where, by slight abuse of notation, we denote by $\langle, \rangle$ the obvious extension of the symplectic form on $V^r$. If $L$ is the dual of an integral lattice, then it is standard that one can find a symmetric positive-definite bilinear form on $V$, which is $\Z$-valued on $L$. In any case, we can define a polarized abelian variety starting from $L$, as above. As the action of $\Gamma$ on $V^r(\R)$ induces isomorphisms between abelian varieties in the same $\Gamma$-orbit, we may identify the fiber $A^r(s):=\varphi^{-1}(s)$, for $s\in S$, with $A^r(x)$ for any $x\in \mathfrak{X}$ above $s$.

As above, we may define complex varieties, for $v\in L^+,$
\[
\mathcal{A}^{(1)}_v:=(\Gamma_v\ltimes L_v)\backslash (\mathfrak{X}_v\times V_v(\R))
\]
and
\[
\mathcal{A}_v:=\mathcal{A}^{(r)}_v:=\underbrace{\mathcal{A}_{v}^{(1)}\times_{S_{v}}\cdots\times_{S_{v}}\mathcal{A}_{v}^{(1)}}_\text{$r$ factors}.
\]

The natural inclusions $\mathfrak{X}_v\hookrightarrow \mathfrak{X}$ and $V_v\hookrightarrow V$ induce a degree one morphism $h_v: \mathcal{A}_v\to \mathcal{A}$ which is compatible with the fiber structure in the sense that the following diagram commutes:
\[
\begin{tikzcd}
A^r_v(s_v)\MySymb{dr} \arrow[d, "j_v"] \arrow[r]     & \mathcal{A}_v \MySymb{dr} \arrow[r,"\varphi_v"] \arrow[d, "h_v"]    & S_v \arrow[d, "\iota_v"]\\
A^r(\iota_v(s_v)) \arrow[r]                         & \mathcal{A} \arrow[r,"\varphi"]                       & S
\end{tikzcd}
\]
If $r=0$ then we take $\mathcal{A}_v=S_v$ and $h_v=\iota_v$.

\subsection{Compactifications}
Let $\overline{S}^{\text{BB}}$ denote the Baily-Borel compactification of $S$. It is possible to find smooth projective toroidal compactifications $\tilde{\mathcal{A}}$ and $\tilde{S}$ of $\mathcal{A}$ and $S$, respectively, so that the following diagram commutes:
\[
\begin{tikzcd}[row sep=large, column sep=large, scale cd=1.3]
\mathcal{A} \MySymb{dr} \arrow[d, "\varphi"] \arrow[r,"\iota"]     & \tilde{\mathcal{A}} \MySymb{dr} \ar[equal]{r} \arrow[d, "\tilde{\varphi}"]    & \tilde{\mathcal{A}} \arrow[d, "\overline{\varphi}"]\\
S \arrow[r]                         & \tilde{S} \arrow[r]                       & \overline{S}^{\text{BB}}
\end{tikzcd}
\]
A smooth toroidal compactification $\tilde{S}$ has been constructed in \cite[\S 2]{cogdell}, which has been denoted $S^{\text{reg}}$ in section \S 2 of this manuscript. A complex variety $\tilde{\mathcal{A}}$ such that the diagram commutes has been constructed in general in \cite[p. 136]{namikawa2}, see also \cite[p. 18]{goren} for the case of Picard modular surfaces.
Analogously, we have a smooth compactification $\tilde{\mathcal{A}_v}$ of $\mathcal{A}_v$ for which the analogous diagram commutes, see \cite[p. 160]{gordon2}:
\[
\begin{tikzcd}[row sep=large, column sep=large, scale cd=1.3]
\mathcal{A}_v \MySymb{dr} \arrow[d, "\varphi_v "] \arrow[r]     & \tilde{\mathcal{A}}_v  \MySymb{dr} \ar[equal]{r} \arrow[d, "\tilde{\varphi}_v "]    & \tilde{\mathcal{A}}_v  \arrow[d, "\overline{\varphi}_v "]\\
S_v  \arrow[r]                         & \tilde{S}_v  \arrow[r]                       & \overline{S}^{\text{BB}}_v 
\end{tikzcd}
\]
\noindent The compatibilities of these constructions are summarized in the following commutative diagram:

\[
    \begin{tikzcd}[row sep=1.75em, column sep = 2.5em]
    \mathcal{A}_v \arrow[rr] \arrow[dr, swap,"h_v"] \arrow[dd,swap, "\varphi_v"] &&
    \tilde{\mathcal{A}_v} \arrow[dd, pos=0.4, "\tilde{\varphi}_v"] \arrow[dr,"\tilde{h}_v"] \\
    & \mathcal{A} \arrow[rr] \arrow[dd, pos=0.35, "\varphi"] &&
    \tilde{\mathcal{A}} \arrow[dd,"\tilde{\varphi}"] \\
    S_v \arrow[rr,] \arrow[dr, "\iota_v"] && \tilde{S_v} \arrow[dr, "\tilde{\iota}_v"] \\
    & S \arrow[rr]&& \tilde{S}
    \end{tikzcd}
\]
cf. \textit{loc. cit.} for more details.

\subsection{Cohomology}
By the degeneration of the Leray spectral sequence for $\varphi:\mathcal{A}\to S$, see \cite[p. 268]{gordon1}, we have
\[
\h^i(\mathcal{A},\Q)\cong \bigoplus_{a+b=i}\h^{a,b}(\mathcal{A},\Q),
\]
where $\h^{a,b}(\mathcal{A},\Q)\cong \h^a(S, R^b\varphi_*\Q)$, for $0\le a\le 4$ and $0\le b\le 6r$. Furthermore, the same decomposition is also valid for cohomology with compact support.
Now, let us fix a base point $x_0\in \mathfrak{X}$ with image $s_0\in S$, and identify $\pi_1(S,s_0)$ with $\Gamma$ by letting $\gamma\in \Gamma$ correspond to the homotopy class of the image in $S$ of a path joining $x_0$ to $\gamma x_0$. Similarly, $\h_1(A^r(x_0),\Z)\cong L^r$, and thus
\begin{equation}
\label{2.3.3}
\h^b(A^r(x_0),\Q)\cong \bigwedge^b (V^{\vee})^r,
\end{equation}
where $V^\vee$ is the dual of $V$. With these identifications, the monodromy action of $\pi_1(S,s_0)$ on $\h^b(A^r(x_0),\Q)$ becomes the $\bigwedge^b \sigma^r$-action of $\Gamma$ on $\bigwedge ^b (V^{\vee})^r$. In particular, we have
\begin{equation}
\label{2.3.4}
\h^{a,b}(\mathcal{A},\Q)\cong \h^a(S,\bigwedge ^b (V^{\vee})^r),
\end{equation}
and
\[
\h^{a,b}_c(\mathcal{A},\Q)\cong \h^a_c(\Gamma,\bigwedge ^b (V^{\vee})^r),
\]
since both are isomorphic to $\h^a(S, R^b\varphi_*\Q)$ (resp. $\h^a_c(S, R^b\varphi_*\Q)$), see \cite[p. 267]{gordon1}.
 Thus each $\h^{a,b}(\mathcal{A},\Q)$ may be decomposed according as $\bigwedge ^b (V^{\vee})^r$ decomposes as a $\Gamma$-module. Note that since $\Gamma$ is Zariski-dense in $G$, $\bigwedge ^b (V^{\vee})^r$ decomposes identically under the action of $\Gamma$ or of $G$.
Write 
\[
\bigwedge\nolimits^* (V^\vee)^r:=\bigoplus_{\xi\in \Sigma_r} V_\xi.
\]

\begin{Lemma}[\text{\cite[Lemma 2.1.4]{gordon2}}]
The decomposition
\[
\h^*(\mathcal{A},\Q)\cong \bigoplus_{\substack{0\le i \le 4\\ \xi\in \Sigma_r}}\h^{i}(S,\mathcal{V}_\xi)
\]
is algebraically defined, or motivic, in the sense that for each term on the right hand side there is an idempotent $\mathcal{P}_\xi$ in the $\Q$-algebra of algebraic correspondences on $\mathcal{A}$ (modulo homological equivalence) which induces a projection from $\h^*(\mathcal{A},\Q)$ onto a subspace isomorphic to that term. The same decomposition is also valid for cohomology with compact support.
\end{Lemma}

\begin{Proposition}
For each $r\geq 0$ there is among all the absolutely irreducible constituents of $\bigwedge^{2r} (V^{\vee})^r$ a unique one $(\rho_r, E_{r})$ of maximal dimension. This representation occurs in $\bigwedge^{2r} (V^{\vee})^r$ with multiplicity one, is therefore defined over $\Q$, and it is equivalent to the irreducible representation of highest weight $r$ of $\GU(2,1)$.
\end{Proposition}

\begin{proof}
The proof is the same as that in \cite[p. 16]{gordon}.
\end{proof}

Now let $\h^{2r}(\mathcal{E}_{r}(x),\Q)$ denote the subspace of $\h^{2r}(A^r(x),\Q)$ isomorphic to $E_{r}$ via (\ref{2.3.3}), and let $\h^{2r+2}(\mathcal{M}, \Q)$ denote the subspace of $\h^{2,2r}(\mathcal{A},\Q)$ isomorphic to $\h^2(S, E_{r})$ via (\ref{2.3.4}).

\begin{Proposition}
The representation-theoretic component of $\h^{2r+2}(\mathcal{A},\Q)$ corresponding to $\h^2(S, E_{r})$ is an algebraically defined subspace of $\h^{2r+2}(\mathcal{A},\Q)$.
\end{Proposition}

\begin{proof}
This follows from the above lemma.
\end{proof}

\begin{Theorem}[\text{\cite[Proposition 2.5]{gordon2}}]
Let $$S^\infty:=\overline{S}^{\emph{BB}}-S$$
be the cuspidal divisor of $\overline{S}^{\emph{BB}}$. Then
\[
\h^{*}(\tilde{\mathcal{A}},\Q)\cong \bigoplus_{\xi\in \Sigma_r}\emph{IH}^i(\overline{S}^{\emph{BB}}, \mathcal{V}_\xi)\oplus \bigoplus_\alpha \h^0(S^\infty,\mathcal{U}_\alpha)[j_\alpha],
\]
for some sheaves $\mathcal{U}_\alpha$ on $S^\infty$ and integers $j_\alpha$. Furthermore the splitting between cohomology supported on $\overline{S}^{\emph{BB}}$ and cohomology supported at the cusps can be taken to be orthogonal with respect to Poincaré duality.
\end{Theorem}

Here, $\text{IH}^*$ denotes intersection cohomology, which in our case reduces to the definition in \cite[Equation (2.3.1)]{gordon2}. We see this result as the analogue in higher weight of the splitting in homology of the smooth compactified Picard modular surface, \cite[p. 125]{cogdell}.

\subsection{Higher weight cycles}
%For $v\in L$ with $(v,v)>0$, let $$\mathfrak{n}_L(v)=\gcd\{\text{tr}_{K/\Q}(u,v): u\in L\}\in \Z.$$
In the following, we will fix $r=2k$, for $k$ a positive integer.
Let $E_{k,k}$ be the irreducible representation of $\SU(2,1)$ obtained from $E_{2k}$ by the map
\[
\bigwedge\nolimits^{4k}(V\oplus V)^k\to \bigwedge\nolimits^{4k}(V\oplus V^\vee)^k
\]
induced by the polarization of $\mathcal{A}$.

The aim of this subsection is to define higher weight cycles in the cohomology of Kuga-Sato varieties attached to $S$. These are cohomology classes
\[
C^k_v\in \h^{4k+2}(\mathcal{W}_k,\Q)
\] 
for 
\[
\mathcal{W}_k:=\widetilde{\mathcal{A}^{(k)}\times_{S} (\mathcal{A}^{(k)})^\vee}\to \tilde{S},
\]
with the following properties:
\begin{itemize}
\item They project to middle-degree classes in $\h^2(S,E_{k,k})$ via a unique projector $\mathcal{P}_k$ in the ring of correspondences of $\mathcal{W}_k$.
\item They are orthogonal to the boundary divisor $S^\infty$, and hence project to $\text{IH}^2(\overline{S}^{\emph{BB}},E_{k,k})$ via $\mathcal{P}_k$.
\item They are supported on the special cycles $C_v^c$.
\end{itemize}

To make the following definitions, we apply the results obtained by Gordon in \cite{gordon2}, recalled in the previous subsection.
\begin{definition}
\begin{enumerate}[(i)]
\item For $k>0$, $v\in L^+$ and $x\in \mathfrak{X}_v$, let
\[
C_v^k(x):=\mathcal{P}_k(x)\circ (j_v)_*(\mathcal{W}_k),
\]
be the algebraic cycle of higher weight on $W_k(x):=A^{(k)}\times (A^{(k)})^\vee(x)$ which represents the projection to $\h^{4k}(\mathcal{E}_{k,k}(x),\Q)$ of the image in $\h^{4k}(W_k(x),\Q)$ of the fundamental class of $A_v^k(x)\times (A_v^k)^\vee(x)$.
\item For each $k>0$ and $v\in L^+$,
\[
C_v^k:=\mathcal{P}_k\circ\iota\circ\circ h_v(\mathcal{A}_v^{(k)}\times_{S_v}(\mathcal{A}_v^{(k)})^\vee)
\]
is an algebraic cycle of higher weight on $\mathcal{W}_k$ which represents a class in $\text{IH}^{4k+2}(\mathcal{M},\Q)$.
\item For each positive integer $n$, we define composite cycles of higher weight 
\[
C_n^k:=\sum_{v\in \Gamma\backslash L(n)} d(L_v)^{-k}d(L_v^\vee)^{-k}C_v^k,
\]
where $d(L_v)$ and $d(L_v^\vee)$ are the discriminants of the lattices $L_v:=L\cap V_v$ and of its dual, respectively.
\item Formally, we also let
\[
C_0^k:=\begin{cases}
\frac{1}{2}c_1(S)\;\; &\text{if } k=0,\\
0\;\; &\text{otherwise},
\end{cases}
\]
where $c_1(S)$ is the first Chern class of $S$.
\item When $k=0$, set $C^0_v:=C_v^c$ and $C^0_n:=C_n^c$.
\end{enumerate}
\end{definition}

Next, we wish to prove modularity of these higher weight cycles, and for this we will compare them to the cycles with coefficients defined in \cite[p. 51]{little}.

As in \cite[p. 23]{gordon}, we shall write $C_v^k(x)$ in coordinates, exploiting an explicit factorization of the projection $\mathcal{P}_k(x)$.
\begin{Lemma}
We have a factorization
\[
\mathcal{P}_k:\bigwedge\nolimits^{4k}(V\oplus V^\vee)^k\twoheadrightarrow V^{\otimes k}\otimes (V^\vee)^{\otimes k}\twoheadrightarrow E_{k,k}.
\]
\end{Lemma}
\begin{proof}
This is obvious from the definition of $E_{k,k}$ as a quotient of $V^{\otimes k}\otimes (V^\vee)^{\otimes k}$, see \cite[p. 40]{little}, and the definition of the above map $\mathcal{P}_k$ as the natural projection onto a subspace, see \cite[p. 267]{gordon1}.
\end{proof}
We will mainly be interested in the first projection
$\bigwedge\nolimits^{4k}(V\oplus V^\vee)^k\twoheadrightarrow V^{\otimes k}\otimes (V^\vee)^{\otimes k}$. Namely, we want to project via $$\bigwedge\nolimits^{4k}(V\oplus V^\vee)^k\cong \h^{4k}(W_k(x),\Q)\twoheadrightarrow V^{\otimes k}\otimes (V^\vee)^{k}$$ the push-forward $j_{v,*}(\eta_v)$ of the fundamental class of $A_v^k(x)\times (A_v^k)^\vee(x)$ induced by the inclusion
\[
j_v:A_v^k(x)\times (A_v^k)^\vee(x)\hookrightarrow W_k(x).
\]
\begin{Lemma} For $v\in L^+$, write $v_1^i,v_2^i$ for its components in the real symplectic space $V^i:=V\;\; (i=1,\dots, k)$, then
\[
j_{v,*}(\eta_v)=d(L_v)^{k}d(L_v^\vee)^{k}\bigwedge_{i=1}^{k}v^i_1\wedge v_2^i\wedge v_1^{i,*}\wedge v_2^{i,*}\in \h^{4k}(W_k(x),\Q).
\]
\end{Lemma}

\begin{proof}
By definition of the Gysin morphism $j_{v,*}$, we have a commutative diagram:
\[
\begin{tikzcd}[row sep=large, column sep=large]
\h^0(A_v^k(x)\times (A_v^k)^\vee(x),\Q) \arrow[d, "\cong"] \arrow[r,"j_{v,*}"]    &   \h^{4k}(W_k(x),\Q)  \\
\h_{8k}(A_v^k(x)\times (A_v^k)^\vee(x),\Q) \arrow[r,"j_{v,*}"] & \h_{8k}(W_k(x),\Q) \arrow[u, "\cong"].
\end{tikzcd}
\]
If we write $V_v^i=\langle e_1^i,e_2^i,e_3^i,e_4^i\rangle_\R$, then the image of $\eta_v$ in $$\h_{8k}(A_v^k(x)\times (A_v^k)^\vee(x),\Q)\cong \bigwedge^{8k} (L_v\oplus L_v^\vee)^k$$ is
\[
d(L_v)^{k}d(L_v^\vee)^{k}\bigwedge_{i=1}^{k}\bigwedge_{l=1}^{4}e^i_l\wedge e_l^{i,*}.
\]
Via the inclusion $V_v^i\subset V^i$, we identify the above with its image in $\h_{8k}(W_k(x),\Q)$. From this the result follows, since $V^i=V^i_v\oplus \langle v_1^i,v_2^i\rangle_\R$.
\end{proof}
By the above lemma, it is clear that, by identifying $v_1^i\wedge v_2^i$ with $v^i$, the image of $j_{v,*}(\eta_v)$ in $ V^{\otimes k}\otimes (V^\vee)^{\otimes k}$ will be
\[
d(L_v)^{k}d(L_v^\vee)^{k}v^{\otimes k}\otimes (v^*)^{\otimes k}.
\]
For each $v\in L^+$, let $C_{v,[k,k]}\in \h^2(S,E_{k,k})$ be the cycle with coefficients defined in \cite[p. 51]{little}. We give a brief overview of the construction in \textit{loc. cit.}.

Given a $G$-representation $E$, and any positive $v\in V$, consider the bundle locally given by the projection:
\[
C_v\times_{\Gamma_v} E\to C_v.
\]
For any $\Gamma_v$-invariant vector $w$ in $E$, we may write sections $s_w$ of the bundle as
\[
s_w(z)=(z,w).
\]
For simplicity we write $C_v\otimes w:= C_v\otimes s_w$.  Fixing $E=V^l\otimes (V^\vee)^{l'}$, we have the natural choice $w=v^l\otimes (v^*)^{l'}$: this is a constant and thus parallel section. We quote a result from Little's thesis:
\begin{Proposition}[\text{\cite[Proposition 4.1.3]{little}}]
Fix non-negative integers $l,l'$ and a positive vector $v\in V$. We then
define the special cycle with coefficients in $V^l\otimes (V^\vee)^{l'}$ as follows:
\[
C_{v,l,l'}:=C_v\otimes v^l\otimes (v^*)^{l'}.
\]
Similarly, we then define the special cycle with coefficients for the representation $\mathcal{H}^{l,l'}(V)$ (see \emph{\cite[p. 40]{little}}) as
\[
C_{v,[l,l']}:=C_v\otimes \pi_{\mathcal{H}}(v^l\otimes (v^*)^{l'})
\]
These are cycles - namely, they are closed - and so in particular represent classes in
homology:
\[
C_{v,l,l'}\in \h_2(S,\partial S, \widetilde{V^l\otimes (V^\vee)^{l'}})\;\;\;\;\;\;\; C_{v,[l,l']}\in \h_2(S,\partial S, \widetilde{\mathcal{H}^{l,l'}(V)}).
\]
\end{Proposition}
\noindent For the comparison with the cycles on the Kuga-Sato variety, we will be interested with the case $l=l'=k$, for which, essentially by definition of the representation $E_{k,k}$ and by \cite[Theorem 3.1.2]{little}, we have 
\[
\mathcal{H}^{l,l'}(V)=\mathcal{H}^{k,k}(V)\cong E_{k,k}.
\]
Let $(,)$ be the natural pairing on cohomology induced by Poincaré duality, cf. \cite[p. 46]{little}, and let 
\[
\pi_{\mathcal{H}}: V^{\otimes k}\otimes (V^\vee)^{\otimes k}\twoheadrightarrow  \mathcal{H}^{k,k}(V) \cong E_{k,k}
\]
be the natural projection.
\begin{Corollary}
For any compactly supported $\varphi\in \h^2(S,E_{k,k})$, we have the equality 
\[
(C_{v,[k,k]},\varphi)=d(L_v)^{-k}d(L_v^\vee)^{-k}(C_v^k,\varphi).
\]
\end{Corollary}
\begin{proof}
Indeed, by definition of the natural projection onto a subspace
\[
\pi_{\mathcal{H}}:V^{\otimes k}\otimes (V^\vee)^{\otimes k}\twoheadrightarrow E_{k,k},
\]
the section
\[
\pi_{\mathcal{H}}(v^{\otimes k}\otimes (v^*)^{\otimes k})\in E_{k,k}
\]
coincides with the section defined in \cite[p. 51]{little}.
\end{proof}

\begin{Theorem}
\label{main2}
The theta series with coefficients in the cohomology of the Kuga-Sato variety $\mathcal{W}_k$
\[
\sum_{n\geq 0}C_n^k q^n\in \emph{H}^{4k+2}(\mathcal{M},\Q)[[q]]
\]
is modular of weight $2k+3$, i.e., for $\varphi\in \emph{H}^{8k+2}(\mathcal{W}_k,\Q)$ with compact support, the $q$-expansion
\[
\sum_{n\geq 0}(C_n^k, \varphi) q^n
\]
is an elliptic modular form of weight $2k+3$.
\end{Theorem}

\begin{proof}
This follows from the previous corollary and from \cite[Theorem $6.4.1$]{little}.
\end{proof}

\begin{Remark}
We expect the above result to hold for any $\varphi\in \text{H}^{8k+2}(\mathcal{W}_k,\Q)$. Indeed, it seems likely that Little's procedure of capping the cycles with coefficients should algebraically correspond to the canonical projection to intersection cohomology of the cycles on the Kuga-Sato variety.   
\end{Remark}
\noindent It follows that, for $k$ a positive integer, the adjoint Kudla lift in weight $2k+4$ is given by integration on the higher weight cycles $C_v^k$. More precisely, let $$j_{k,k}:\mathscr{S}(2k+4,\Gamma)\to \h^2(S,\mathcal{H}^{k,k}(V))$$ be the holomorphic Eichler-Shimura mapping defined in \cite[p. 246]{hidashimura}, and given explicitly by 
\[
j_{k,k}(F)=F\otimes \underline{\mathfrak{z}}^{\otimes k}\otimes (\underline{\mathfrak{z}}^*)^{\otimes k}d\mathfrak{z}.
\] 
Then we have the following:
\begin{Corollary}
\label{relation}
Let $F\in \mathscr{S}(2k+4,\Gamma)$. Then we have the formulas
\begin{enumerate}[(i)]
\item $$(F,\theta_v)=c_{2k+4}(L)d(L_v)^{-k}d(L_v^\vee)^{-k}\Big[\int_{C_v^k}j_{k,k}(F)\Big],$$
\item 
\[
\mathcal{L}^\dagger(F)=c_{2k+4}(L)\sum_{v\in L^+}d(L_v)^{-k}d(L_v^\vee)^{-k}\Big[\int_{C_v^k}j_{k,k}(F)\Big]q^{(v,v)},
\]
\end{enumerate}

\end{Corollary}
\begin{proof}
This follows from the above and from the explicit formula in Proposition \ref{analogue}, where it is essential to remark that since $\pi_{\mathcal{H}}$ is an orthogonal projection with respect to the pairing in $V$, see \cite[p. 26]{gordon}, we have
\[
(\underline{\mathfrak{z}}^{\otimes k}\otimes (\underline{\mathfrak{z}}^*)^{\otimes k},\pi_{\mathcal{H}}(v^{\otimes k}\otimes (v^*)^{\otimes k}))=(\underline{\mathfrak{z}}^{\otimes k}\otimes (\underline{\mathfrak{z}}^*)^{\otimes k},v^{\otimes k}\otimes (v^*)^{\otimes k})=(\underline{\mathfrak{z}},v)^{2k}.
\]
\end{proof}

\section{A $\Lambda$-adic Family of Special Cycles}
In this section we construct a $p$-adic analytic family of unitary special cycles, by applying Loeffler's methods \cite{loeffler}. Later, we modify these cohomology classes slightly in order to be able to compute their intersection multiplicities under a suitably twisted "big" Poincaré pairing, which can be seen as a generalization of a construction of Ohta \cite[Proposition 4.1.13]{ohta}, see also \cite[Proposition 7.13]{fornea} or \cite[Lemma 1.1]{darmon} for the same construction in a different setting. The $\Lambda$-adic intersection multiplicities are then packaged in a formal $q$-expansion, whose specializations at arithmetic points of weight $k\geq 3$ are proven to be modular, see Theorem \ref{bigcog} below, providing us with a $\Lambda$-adic version of the Cogdell theorem:
\begin{namedthm}{Main Theorem}
Let $$\zeta_\infty\in e\h_{\emph{Iw}}^2(Q^0_G,\Z_p)$$ be a cuspidal big cohomology class such that $\zeta_1$ is orthogonal to the Chern class on $S_{G}(V_1)$, the $\Lambda$-adic $q$-expansion
\[
\Phi_{\zeta}(z):=\sum_{n\geq 1}[\xi_{n,\infty}^{\emph{n.o.}}, \zeta_\infty] q^n\in \Lambda[[q]]
\]
interpolates the modular forms constructed by Cogdell and their higher weight analogues, at each odd weight $k\geq 3$ and level $r> 0$, i.e.
\[
\nu_{k,r}(\Phi_{\zeta}(z))\in S_{k}(\Gamma_1(p^rD)).
\]
\end{namedthm}
Here, the big classes $\xi_{n,\infty}^{\text{n.o.}}$ are constructed \textit{ad hoc} to interpolate the cycles $C_n^k$ for each $k$.
\subsection{Loeffler's machinery}
\label{machinery}
We recall the setting of \cite[\S 4]{loeffler} in order to apply the formalism to our case in the next section.
We suppose to have an inclusion $\iota:H\hookrightarrow G$ of reductive group schemes over $\Z_p$. Let us fix a Borel and a maximal torus $B_G, T_G\subseteq G$ so that their intersection $B_H, T_H$ with $H$ are a Borel and a maximal torus in $H$.

Fix a parabolic $Q_G\supseteq B_G$. By a mirabolic subgroup of $G$ we mean a group of the form
\[
Q^0_G=N_G\cdot L^0_G,
\]
where $N_G$ and $L_G$ are respectively the unipotent and the Levi appearing in the Levi decomposition of $Q_G$, and $L^0_G\unlhd L_G$. We define mirabolic subgroups of $H$ in a similar fashion.

Let $Q^0_H$ be a mirabolic in $H$, and $Q^0_G$ a mirabolic in $G$. We consider the left action of $G$ on the flag variety $\mathcal{F}=G/\overline{Q}_G$, where $\overline{Q}_G$ is the opposite of $Q_G$, relative to our fixed maximal torus $T_G$. We may also let $Q^0_H$ act on $\mathcal{F}$ via the embedding $\iota$. We assume that there exist an element $u\in G(\Z_p)$ with the following properties:
\begin{enumerate}[(A)]
\item The $Q^0_H$-orbit of $u$ is open in $\mathcal{F}$,
\item We have
\[
u^{-1}Q^0_Hu\cap \overline{Q}_G\subseteq \overline{Q}^0_G.
\]
\end{enumerate}

As condition (B) is always satisfied by taking $L^0_G=L_G$, the most interesting results will be obtained by taking a smaller $L^0_G$.

We fix a cocharacter $\eta\in X_\bullet(T_G)$ which factors through $Z(L_G)$ and which is strictly dominant. We set $\tau=\eta(p)\in G(\Z_p)$. Following \cite{loeffler}, we define level subgroups:
\begin{definition} Fix $K^{(p)}\subset G(\hat{\Z}^{(p)})$ a neat open compact subgroup and for $r\geq 0$ let
\begin{align*}
K_r &=K^{(p)}\cdot\{g\in G(\Z_p):\tau^{-r}g\tau^r\in G(\Z_p), (g\!\!\mod p^r)\in \overline{Q}^0_G(\Z/p^r\Z)\},\\
K'_r &=K^{(p)}\cdot\{g\in G(\Z_p):\tau^{-r-1}g\tau^{r+1}\in G(\Z_p), (g\!\!\mod p^r)\in \overline{Q}^0_G(\Z/p^r\Z)\},\\
V_r&=\tau^{-r}K_r\tau^r.
\end{align*}
\end{definition}

Now suppose that $M_G$ is a Cartesian cohomology functor for $G$, and $M_H$ is a Cartesian cohomology functor for $H$, see \cite[\S3]{loeffler} for a precise definition. For our purposes, the case of Betti cohomology will suffice, but we state the next result in all generality:

\begin{Theorem}\emph{\cite[Proposition 4.5.1]{loeffler}}
\label{loeffler}
Under assumptions \emph{(A)} and \emph{(B)}, the following diagram is commutative
\[
\begin{tikzpicture}[baseline=(current bounding box.center), 
    every node/.style={scale=1}, 
    text height=1.5ex, text depth=0.25ex]

% Adjusted spacing
\def\xshift{6} % Horizontal spacing
\def\yshift{1.5}  % Vertical spacing

% Nodes
\node (A1) at (0,2*\yshift) {$M_{H,\mathrm{Iw}}\big(Q_H^0 \cap u K_{r+1} u^{-1} \big)$};
\node (A2) at (\xshift,2*\yshift) {$M_G(K_{r+1})$};
\node (B1) at (0,\yshift) {$M_{H,\mathrm{Iw}}\big(Q_H^0 \cap u K_r' u^{-1} \big)$};
\node (B2) at (\xshift,\yshift) {$M_G(K_r')$};
\node (B3) at (1.6*\xshift,\yshift) {$M_G(\tau^{-1} K_r' \tau)$};
\node (C1) at (0,0) {$M_{H,\mathrm{Iw}}\big(Q_H^0 \cap u K_r u^{-1} \big)$};
\node (C2) at (\xshift,0) {$M_G(K_r)$};
\node (C3) at (1.6*\xshift,0) {$M_G(K_r)$};

% Horizontal arrows
\draw[->] (A1) -- (A2) node[midway, above] {$[u]_{\star}$};
\draw[->] (B1) -- (B2) node[midway, above] {$[u]_{\star}$};
\draw[->] (B2) -- (B3) node[midway, above] {$[\tau]_{\star}$};
\draw[->] (C1) -- (C2) node[midway, above] {$[u]_{\star}$};
\draw[dotted, ->] (C2) -- (C3);

% Vertical arrows
\draw[double equal sign distance] (A1) -- (B1);
\draw[->] (C1) -- (B1);
\draw[->] (A2) -- (B2);
\draw[->] (C2) -- (B2);
\draw[->] (B3) -- (C3);

% Improved curved dashed arrow
\draw[dashed, ->]plot[smooth, tension=.9] coordinates {(\xshift+1,2*\yshift)(\xshift+3,2*\yshift-0.45)(11.5,1.5)(1.6*\xshift+0.5,0.4)};
%\draw[dashed, ->]  to[out=10, in=10, looseness=1.5] (C3.north) 
    \node[circle, minimum size = 4pt, inner sep =0pt]  (c2u) at (9.75, 2.92) {$[\tau]_{K_{r+1}, K_r, \star}$};

\end{tikzpicture}
\]

\end{Theorem}

\subsection{The $\U(2,1)$-case}

In this section we fix $G=\GL_3$ as a group scheme over $\Z_p$ and $H=\GL_2$ a subgroup of $G$ via
\begin{equation*}
g\mapsto \begin{pmatrix}
g &\\
& 1
\end{pmatrix},
\end{equation*}
and take $Q^0_H:=H$ as the mirabolic in $H$, so that we are in the "Hecke-type" case, following the terminology of \cite[\S 5]{loeffler}. For $G$, we set $Q_G=B_G$, the upper triangular Borel, and
\[
Q^0_G=\{g=(g_{ij})\in B_G\mid g_{11}=g_{33}=1\}.
\]
Denote by $N_G$ the upper unipotent of $G$, so that $B_G=T_GN_G$. We also let $\mathcal{H}/\Q$ and $\mathcal{G}/\Q$ be suitable globalizations of $H$ and $G$, respectively, i.e.
\[
\mathcal{H}(\Q_p)\cong H(\Q_p),\;\;\; \mathcal{G}(\Q_p)\cong G(\Q_p).
\]
As this will be our primary interest, we also suppose that $\mathcal{H}(\R)\cong \U(1,1)(\R)$ and $\mathcal{G}(\R)\cong \U(2,1)(\R)$. Switching to unitary groups, from the unitary similitude groups considered in the previous sections, will not matter, as we do not use the PEL structure here, and only restrict our attention to Betti cohomology.  Fix $$u=\begin{pmatrix}
1&0&1\\
&1&1\\
&&1
\end{pmatrix}\in N_G.$$
\begin{Lemma}
\label{46}
\begin{enumerate} 
\item The map 
\begin{align*}
H\times \overline{B}_G&\to G\\
(h,\overline{b})&\mapsto hu\overline{b}
\end{align*}
is an open immersion.
\item $u^{-1}Hu\cap \overline{B}_G\subset \overline{Q}^0_G$.
\end{enumerate}
\end{Lemma}

\begin{proof}
Let us begin by computing $u^{-1}Hu\cap \overline{B}_G$. For $\begin{pmatrix}
a&b\\
c&d
\end{pmatrix}\in H$, we have
\begin{align*}
\begin{pmatrix}
1&0&-1\\
&1&-1\\
&&1
\end{pmatrix} \begin{pmatrix}
a&b&\\
c&d&\\
&&1
\end{pmatrix} \begin{pmatrix}
1&0&1\\
&1&1\\
&&1
\end{pmatrix}=\begin{pmatrix}
a&b&a-1+b\\c&d&c+d-1\\
&&1
\end{pmatrix},
\end{align*}
and for this matrix to belong to $\overline{B}_G$, one must have $b=0, a=1, c=1-d$, so that the stabilizer of $u\overline{B}_G\in G/\overline{B}_G=\mathcal{F}$ is $1$-dimensional and $u^{-1}Hu\cap \overline{B}_G\subset \overline{Q}^0_G$.
Finally, the $H$-orbit of $u$ is open in the Borel flag because its dimension is $\dim H-1$, which equals the dimension of $\mathcal{F}$.

\end{proof}

We fix the matrix $$\tau=\begin{pmatrix}
p^2&&\\
&p&\\
&&1
\end{pmatrix}\in \GL_3(\Z_p),$$
and, for this element, we have corresponding level subgroups as in \S \ref{machinery}. For $a,b,c,d,e,f,g,h,i\in \Z_p$, we see
\[
\tau^{-r}\begin{pmatrix}
a&b&c\\
d&e&f\\
g&h&i
\end{pmatrix}\tau^{r}=\tau^{-r}\begin{pmatrix}
p^{2r}a&p^{r}b&c\\
p^{2r}d&p^{r}e&f\\
p^{2r}g&p^{r}h&i
\end{pmatrix}=\begin{pmatrix}
a&p^{-r}b&p^{-2r}c\\
p^{r}d&e&p^{-r}f\\
p^{2r}g&p^{r}h&i
\end{pmatrix},
\]
which shows that 
\[
\begin{pmatrix}
a&b&c\\
d&e&f\\
g&h&i
\end{pmatrix}\in (K_r)_p\Leftrightarrow a,i\equiv 1\!\!\mod p^r,b\equiv f\equiv 0\!\!\mod p^r, c\equiv 0\!\!\mod p^{2r}.
\]

Furthermore, we have inclusions:
\[
K_{r+1}\subset K_r'=K_r\cap \tau K_r\tau^{-1}\subset K_r,
\]
which imply the existence of a morphism
\[
[\tau]_{r+1}: S_G(K_{r+1})\xrightarrow{\text{pr}} S_G(K'_r)\xrightarrow{\tau}S_G(\tau^{-1}K'_r\tau)=S_G(\tau^{-1}K_r\tau\cap K_r)\xrightarrow{\text{pr}}S_G(K_r),
\]
where $S_G(\cdot)$ denotes the Picard modular surface for an unspecified level, and $\text{pr}$ denote natural projections. Note that, by construction, the normalized Hecke correspondence associated to the double coset $[K_r\tau^{-1}K_r]$ is $U_p^*$, the adjoint of $U_p$, which is associated to $[K_r\tau K_r]$, see \cite[p. 4]{rockwood}.

As in \cite[p. 10]{loeffler}, having fixed a choice of globalization $\mathcal{H}/\Q$ of $H$ such that $\mathcal{H}(\R)\cong \U(1,1)(\R)$ and an embedding $\mathcal{H}\hookrightarrow \mathcal{G}$, we take the identity class $z_H\in \h^0_{\text{Iw}}(H,\Z_p)$, and define $z_{G,r}$ to be the image of $z_H$ under the map
\[
\h^0_{\text{Iw}}(H,\Z_p)\xrightarrow{\text{pr}^*}\h^0_{\text{Iw}}(H\cap u^{-1}K_ru,\Z_p)\xrightarrow{[u]_*} \h^2_{\text{Iw}}(K_r,\Z_p).
\]

\begin{Proposition}\emph{\cite[p. 10]{loeffler}}
We have $[\tau]_{r+1,*}(z_{G,r+1})=U_p^*\cdot z_{G,r}$.
\end{Proposition}
\begin{proof}
This is a simple consequence of Theorem \ref{loeffler}, where assumptions (A) and (B) hold because of Lemma \ref{46}.
\end{proof}

We will find it convenient to work with the twisted classes
\[
\xi_{G,r}^\sharp:=[\tau^r]_* z_{G,r}\in \h^2_{\text{Iw}}(V_r,\Z_p),
\]
which satisfy the analogous compatibility:
\begin{Proposition}\emph{\cite[p. 11]{loeffler}}
We have $\emph{pr}_{r+1,*}(\xi_{G,r+1}^\sharp)=U_p^*\cdot  \xi_{G,r}^\sharp$.
\end{Proposition}

\begin{Remark}
We use here the normalization introduced in \cite[p. 3]{rockwood}, so that all correspondences are integrally defined.
\end{Remark}

The following lemma, will allow us to construct an Atkin-Lehner operator:
\begin{Lemma}
\begin{enumerate}[(i)]
\item The matrix \[
\gamma_r:=\begin{pmatrix}
 & & p^{2r}\\
 & p^{r}&\\
 1 & &  
\end{pmatrix}
\] normalizes $K_r$,
\item The transpose matrix ${^t}\gamma_r=\tau^{-r}\gamma_r\tau^{r}$ normalizes $V_r$,
\item $\gamma_r^2={^t}\gamma^{2}_r=p^{2r}\mathbb{I}_{3}$,
\item $\gamma_r\tau\gamma_r^{-1}={^t}\gamma_r\tau{^t}\gamma_r^{-1}=\begin{pmatrix}
1 & &\\
& p &\\
& &p^{2}
\end{pmatrix}=p^{2}\tau^{-1}$.
\end{enumerate}
\end{Lemma}

\begin{proof}
    We just remark that $\gamma_r=\tau^rw_0$, where $w_0$ is the long Weyl element, so that 
    \begin{align*}
        \gamma_r^{-1}K_r\gamma_r=w_0\tau^{-r}K_r\tau^rw_0=w_0 {^t} K_rw_0=K_r,
    \end{align*}
    where we have used that $\tau^{-r}K_r\tau^r=V_r= {^t}K_r$. The rest of the lemma easily follows by similar arguments. 
\end{proof}

\noindent Let $S_G$ be the Picard modular surface associated to the globalization $\mathcal{G}/\Q=\U(2,1)/\Q$ of the group $G$.
As a consequence of the above lemma, the associated morphisms of varieties
\[
\lambda_r: S_G(V_r)\to S_G(V_r)
\]
satisfy
\[
U_p=\lambda_r^*\cdot U_p^*\cdot \lambda_{r,*},
\]
and we thus see $\lambda_r$ as a generalization of the classical Atkin-Lehner operator.
Thus, we may twist our cohomology classes by the Atkin-Lehner correspondence to get 
\[
\xi_{G,r}:=\lambda_r^*\xi_{G,r}^\sharp\in \h^2_{}(V_r,\Z_p),
\]
and we define projections $\pi_{r+1}:S_{G}(V_{r+1})\to S_G(V_{r})$ by the formulas
\[
\pi_{r+1}\circ \lambda_{r+1}^{-1}=\lambda_r^{-1}\circ \text{pr}_{r+1},\,\,\,\,\,\,\,\,\, \text{pr}_{r+1}\circ \lambda_{r+1}^{-1}=\lambda_r^{-1}\circ \pi_{r+1}.
\]
\begin{Proposition}
\label{compatibility}
We have
\[
\pi_{r+1,*}(\xi_{G,r+1})=U_p\cdot  \xi_{G,r}.
\]
\end{Proposition}

\begin{proof}
This follows easily from the identity $U_p=\lambda_r^*\cdot U_p^*\cdot \lambda_{r,*}$, so that 
\begin{align*}
\pi_{r+1,*}(\xi_{G,r+1})&=\pi_{r+1,*}(\lambda_{r+1}^*(\xi_{G,r+1}^\sharp))\\
&=\lambda_r^*(\text{pr}_{r+1,*}(\xi_{G,r+1}^\sharp))=\lambda_r^*(U_p^*\xi_{G,r}^\sharp)\\
&=U_p\lambda_{r}^*  \xi_{G,r}^\sharp=U_p\cdot  \xi_{G,r}.
\end{align*}
\end{proof}

\begin{Remark}
The projections $\pi_{r+1}$ and $\text{pr}_{r+1}$ are the two natural degeneracy maps $S_{G}(V_{r+1})\to S_G(V_{r})$, see \cite[p. 44]{fornea}.
\end{Remark}
\subsection{Big cohomology classes}
\label{bgc}

By projecting the cycles $\xi_{G,r}$ to the nearly ordinary part of cohomology
\[
\xi_{G,r}^{\text{n.o.}}:=e\xi_{G,r}\in e\h^2(V_r,\Z_p)
\]
with respect to the projector
\[
e:=\lim_{n\to \infty} (U_p)^{n!},
\]
we may finally construct a big cycle
\[
\xi_{G,\infty}^{\text{n.o.}}:=\lim_{_{\substack{\leftarrow\\r}}}\,(U_p)^{-r}\!\cdot\xi_{G,r}^{\text{n.o.}}\in e\h_{\text{Iw}}^2(Q_G^0,\Z_p),
\]
where the inverse limit is taken with respect to the degeneracy maps $\pi_{r+1,*}$.
%Note that the notation on cohomology as above makes sense because of Theorem \ref{main} and Hida's classical control theorem for $\GL_2$.

\subsubsection{The big Pairing}
We define new congruence groups
\begin{align*}
    &K_r^0=\{g=(a_{ij})\in G(\hat{\Z})\mid (a_{ij})_p\equiv 0 \mod p^{(j-i)r} \text{ for } i< j\},\\ &V^0_r=\tau^{-r}K_r\tau^r,
\end{align*}
so that the finite groups $V^0_r/V_r\cong K^0_r/K_r$ act as diamond operators on $S_{G}(K_r)$ and $S_{G}(V_r)$.
By means of the inclusion
\begin{align*}
    \Gamma_r\ni z\mapsto [z\mathbb{I}_{3}]\in V^0_r/V_r,
\end{align*}
we obtain an embedding $$\Gamma_r:=1+p(\Z/p^r\Z)\le (\Z/p^r\Z)^\times\le V^0_r/V_r.$$
We are in the position of defining a $\Lambda_r:=\Z_p[\Gamma_r]$-valued pairing 
\[
\h^2(V_r,\Z_p)\times e\h^2(V_r,\Z_p)\to \Lambda_r,
\]
by the formula:
\[
[x_r,y_r]_r:=\sum_{\sigma\in \Gamma_r}(x_r^\sigma,\lambda_{r,*}U_p^{r}\cdot y_r)_r[\sigma^{-1}],
\]
where $(\cdot,\cdot)_r$ is the usual Poincaré duality on Betti cohomology, and the subscript $c$ on  means that we are considering the submodule of $\h^2(V_r,\Z_p)$ generated by compact special cycles, i.e. those corresponding to Shimura curves.

\begin{Proposition}
\label{pairing}
\begin{enumerate}[(i)]
\item The pairing $[\cdot,\cdot]_r$ is $\Lambda_r$-bilinear and all the Hecke operators away from the level $K^{(p)}$ are self-adjoint with respect to it.
\item Let $p_{r+1}:\Lambda_{r+1}\to \Lambda_r$ be the homomorphism induced by the natural projection $\Gamma_{r+1}\to \Gamma_r$. The diagram
\[
\begin{tikzcd}[row sep=huge, column sep=huge]
\h^2(V_{r+1},\Z_p)\times e\h^2(V_{r+1},\Z_p)_{c} \arrow[d, "\pi_{r+1,*}\times \pi_{r+1,*}"] \arrow[r,""]    &   \Lambda_{r+1} \arrow[d,"p_{r+1}"]  \\
\h^2(V_r,\Z_p)\times e\h^2(V_r,\Z_p)_{c} \arrow[r,""] & \Lambda_r 
\end{tikzcd}
\]
%\[
%\begin{tikzcd}[row sep=huge, column sep=huge]
%a\arrow[d, "\pi_{r+1,*}\times \pi_{r+1,*}"] \arrow[r,"[,]_{r+1}"]    &   \Lambda_{r+1}\arrow[d, "p_{r+1}"]  \\
%a \arrow[r,"[,]_r"] & \Lambda_r
%\end{tikzcd}
%\]
commutes.
\end{enumerate}
\end{Proposition}

\begin{proof}
The proof is completely formal and follows by the same arguments as in \cite[p. 44-45]{fornea}, see also \cite[p. 12]{darmon}.
\end{proof}

The above result allows us to package the pairings $[\cdot,\cdot]_{r}$ into a big pairing valued in $\Lambda=\Z_p[[\Gamma]]$, for $\Gamma=1+p\Z_p$:
\[
\h_{\text{Iw}}^2(Q_G^0,\Z_p)\times e\h_{\text{Iw}}^2(Q_G^0,\Z_p)_{c}\xrightarrow{[\cdot,\cdot]} \Lambda.
\]

\subsection{The $p$-adic variation of Cogdell's theorem}
In this section we obtain a $\Lambda$-adic version of Cogdell's main theorem in \cite{cogdell}. In order to achieve this, we use the families of cycles we constructed in \S\ref{bgc}, and pair them via the big pairing. 

The input of Loeffler's machinery are the compatible identity classes $z_{H,r}\in \h^0(H\cap u^{-1}K_ru,\Z_p)$. As divisors on Picard modular surfaces $S_G(K_r)$, these depend on the choice of a globalization $\mathcal{H}$ of the group $H$, and of its embedding in $\mathcal{G}$. The possible choices for $\mathcal{H}$ such that $\mathcal{H}(\R)\cong\U(1,1)(\R)$ correspond to the choice, for each $r\geq 0$, of positive vectors $v_r$ in fixed lattices $L_r\subset V$ preserved by the level $K_r$. This, on the other hand, is the input for the construction of the special cycles $C_{v_r}$, see \S $2$. If it is clear from the context, we will drop the subscript $r$ from the notation $v_r$, thus only writing $v$.

For each positive vector $v\in L_r$, we denote by $$z_H^v\in \h^0(H\cap u^{-1}K_ru,\Z_p)$$ the class corresponding to the vector $v\in L^+_r$, and for a fixed $r$, we may thus consider the cohomology class
\[
\xi_{G,v,r}^{\text{n.o.}}\in e\h^2(V_{r},\Z_p),
\]
constructed by the input $z_H^v$. For ease of notation, we will simply denote these as $\xi_{v,r}^{\text{n.o.}}$. 

We now assume to have, for each $r\geq 0$, an integral lattice $L_r$ preserved by $K_r$, together with its dual 
\[
L_r^\vee:=\{u\in V\mid \text{tr}(u,w)\in \Z,\,\forall w\in L_r\}.
\]
Then $L_r^\vee/L_r$ is finite, and we furthermore assume that there exists, for each $r$, an element $h_r\in L_r^\vee/L_r$ fixed by $K_r$, and such that $h_{r+1}\equiv h_r\mod L_r$. Note that it is always possible to find such lattices and cosets since we have $K_r\subset K_r^0$ and $K_r^0$ is of the form considered in \cite[\S$4.2$]{iudica}, where such choices are made explicit. We then define, for $h_r$ as above,
\[
L'_r:=\{v\in L_r^\vee\mid v\equiv h_r\mod L_r\}.
\]

Finally, for fixed integers $r\geq 0$ and $n\geq 1$, we define cycles
\[
\xi_{n,r}^{\text{n.o.}}:= \sum_{v\in K_r\backslash L'_r(n)}\xi_{v,r}^{\text{n.o.}},
\]
and 
\[
\xi_{n,\infty}^{\text{n.o.}}:= \lim_{_{\substack{\leftarrow\\r}}}U_p^{-r}\cdot\xi_{n,r}^{\text{n.o.}}\in e\h_{\text{Iw}}^2(Q^0_G,\Z_p).
\]
%\begin{Lemma}
%There exists a $r_0\geq 0$ such that for $r\geq r_0$, the classes $\tilde{\xi}_{n,r}^{\emph{n.o.}}\in e'\h_{\et}^2(V_{r}^1,\Z_p(1))_{\theta}$ are divisible by $[V_{r}:V_{r_0}]$.
%\end{Lemma}
%\begin{proof}
%The $\Z_p$-modules $e'\h_{\et}^2(V_{r}^1,\Z_p(1))_{\theta}$ are well-controlled in $r$ because of Theorem \ref{main} and by Hida's classical horizontal control theorem for $\GL_2$, \cite[Theorem 4.2.37]{hidashimura}, they are free and finitely generated. However, as $r$ increases we have $|K_r\backslash L'_r(n)|\to \infty$. This implies that for $r$ big enough and $v_1,v_2\in L'_r(n)$ which are $K_{r-1}$-equivalent but not $K_r$-equivalent, the cohomology classes $\xi_{v_1,r}^{\text{n.o.}}$ and $\xi_{v_2,r}^{\text{n.o.}}$ coincide. The result follows.
%\end{proof}

%We finally define
%\[
%\xi_{n,r}^{\text{n.o.}}:= [V_{r}:V_{r_0}]^{-1}\sum_{v\in V_r\backslash L_r(n)}\xi_{v,r}^{\text{n.o.}}\in e'\h_{\et}^2(V_{r}^1,\Z_p(1))_{\theta}
%\]

\begin{Theorem}
\label{bigcog}
Given a cuspidal big cohomology class $$\zeta_\infty^{\emph{n.o.}}\in e\h_{\emph{Iw}}^2(Q^0_G,\Z_p)_{c}$$ such that $\zeta^{\emph{n.o.}}_1$ is orthogonal to the Chern class on $S_{G}(V_1)$. Then $\Lambda$-adic $q$-expansion
\[
\Phi_{\zeta_\infty^{\emph{n.o.}}}(z):=\sum_{n\geq 1}[\xi_{n,\infty}^{\emph{n.o.}}, \zeta_\infty^{\emph{n.o.}}] q^n\in \Lambda[[q]]
\]
interpolates the modular forms constructed by Cogdell and the modular forms in \S 5, i.e. for each weight $k\geq 0$ and level $r> 0$, 
\[
\nu_{2k,r}(\Phi_{\zeta_\infty^{\emph{n.o.}}}(z))\in S_{2k+3}(\Gamma_1(p^rD)).
\]
\end{Theorem}
\begin{remark}
In \cite{cogdell} the author considers intersection multiplicities of cycles "capped" at infinity, defined on a smooth compactification of the Picard modular surface. As it is not clear how to extend Theorem \ref{loeffler} to compactifications, in the above theorem we have restricted to the case where the special cycle $\zeta_\infty^{\text{n.o.}}$ is cuspidal, so that any contribution at infinity is in fact killed.
\end{remark}
\begin{proof}
First suppose that $k=0$.
By Proposition \ref{pairing}, the specializations of the coefficients of $\Phi_{\zeta_\infty^{\text{n.o.}}}$ correspond to the pairing at finite levels of the specializations of the big cycles. That is,
\begin{align*}
\nu_{r}([\xi_{n,\infty}^{\text{n.o.}}, \zeta_\infty^{\text{n.o.}}])&=[U^{-r}_p\xi_{n,r}^{\text{n.o.}},U^{-r}_p\zeta_r^{\text{n.o.}}]_r=[\xi_{n,r}^{\text{n.o.}},U^{-2r}_p\zeta_r^{\text{n.o.}}]_r\\
&=[\xi_{n,r},U^{-2r}_p\zeta_r^{\text{n.o.}}]_r=\sum_{\sigma\in \Gamma_r}(\xi_{n,r},(\lambda_{r,*}U_p^{-r}\cdot \zeta_r^{\text{n.o.}})^\sigma)_r[\sigma].
\end{align*}
Now, $\xi_{n,r}$ is a twist of the class 
\[
\sum_{v\in K_r\backslash L'_r(n)} z_{G,r}^v
\]
by the Atkin-Lehner correspondence $\lambda_{r,*}$ and the correspondence $[\tau^r]_*$. As the $z_{G,r}^v$ are Cogdell's special cycles, to conclude we need to apply on both arguments of the pairing the inverse correspondences $\lambda_r^*$ and $[\tau^r]^*$, and we get the result for $k=0$. Now take $k\geq 0$. By \cite[Corollary 3.4 and Theorem 3.3]{hida89}, see also \cite[Theorem 7.8]{nicole}, the specialization in weight $2k$ of the big pairing factors through the moment map $\text{mom}^{k,k}$. Hence, we are left to consider, for each $n>0$, the class
\[
\text{mom}^{k,k}_r(\xi_{n,\infty}^{\text{n.o.}})\in e\h^2(V_r,\mathcal{H}^{k,k}).
\]
By \cite[Theorem 5.2.3]{rockwood} the above class corresponds to the ordinary part of the $n$-th higher weight special cycle; indeed, in our case the branching polynomials correspond, up to normalization, to the projection to the highest weight representation of the powers $v^{k}$ of the positive vectors $v\in L'_r(n)$, as the latter is tautologically invariant under the action of any mirababolic of $H_v=\text{Stab}_G(v)$. As the big pairing is Hecke equivariant, we may proceed as above for the case $k=0$, thus getting the result for positive $k$.
\end{proof}

\section*{Acknowledgments}
%\addcontentsline{toc}{section}{Acknowledgements}

I am extremely grateful to Professor Marc-Hubert Nicole for his great help and support during the writing of this paper. I also wish to thank Maria Rosaria Pati for numerous discussions and Matteo Longo for his helpful comments. %and finally the anonymous referee for the many suggestions which helped me improve the accuracy and legibility of this paper.

\end{document}